\documentclass[10pt]{article}
\usepackage[dvips]{color}

\newlength{\minitwocolumn}
\font\teneufm=eufm10
\font\seveneufm=eufm7
\font\fiveeufm=eufm5
\newfam\eufmfam
\textfont\eufmfam=\teneufm
\scriptfont\eufmfam=\seveneufm
\scriptscriptfont\eufmfam=\fiveeufm

\makeatletter
\@addtoreset{equation}{section}
\makeatother

\newtheorem{thm}{Theorem}[section]
\newtheorem{prop}[thm]{Proposition}

\newtheorem{cor}[thm]{Corollary}

\title{\Large{\bf A BOSONIZATION OF $U_q(\widehat{sl}(M|N))$}}
\begin{document}
\maketitle

\begin{center}
{TAKEO KOJIMA}
\\~\\
{\it
Department of Mathematics and Physics,
Faculty of Engineering, Yamagata University,\\
Jonan 4-3-16, Yonezawa 992-8510, JAPAN}
\\kojima@yz.yamagata-u.ac.jp
\end{center}


\begin{abstract}
A bosonization of the quantum affine superalgebra $U_q(\widehat{sl}(M|N))$ is presented for an arbitrary level $k \in {\bf C}$.
Screening operators that commute with $U_q(\widehat{sl}(M|N))$ are presented for the level $k \neq -M+N$.
\end{abstract}


\section{Introduction}
Bosonization is a powerful method to study representation theory 
of infinite-dimensional algebras \cite{Feigin-Fuchs}
and its
application to 
mathematical physics, such as calculation of correlation functions of exactly solvable models \cite{JM}.
In this article we give a bosonization of the quantum affine superalgebra $U_q(\widehat{sl}(M|N))$ $(M,N=1,2,3,\cdots)$ for an arbitrary level $k \in {\bf C}$,
and give screening operators that commute with $U_q(\widehat{sl}(M|N))$ for the level $k \neq -M+N$.
For level $k=1$, bosonization has been constructed for the quantum affine algebra $U_q(g)$ in many cases $g=(ADE)^{(r)}, (BC)^{(1)}, G_2^{(1)}, \widehat{sl}(M|N),
{osp}(2|2)^{(2)}$ \cite{FJ, B, JKM, Jing1, Jing2, KSU, Z, YZ1, YZ2}.
Bosonization of an arbitrary level $k \in {\bf C}$ is completely different from those of level $k=1$.
For level $k \in {\bf C}$,
bosonization has been constructed only for $U_q(\widehat{sl}(N))$ and $U_q(\widehat{sl}(N|1))$
\cite{W, Feigin-Frenkel, IK, M, S, Kon, AOS1, AOS2, Koj1}.  
In this article we give a higher-rank generalization of the previous works
for the quantum affine superalgebra $U_q(\widehat{sl}(N|1))$ 
including the construction of screening operators
\cite{AOS2, Koj1, ZG, Koj2, Koj3}. 
Representation theory
of the superalgebra is
 much more complicated than non-superalgebra and has 
rich structures \cite{Kac1, Kac2, Kac-Wakimoto, Kac3, Yamane, Palev-Tolstoy}.

The text is organized as follows. 
In Section \ref{Section:2} we recall the Chevalley generators and the Drinfeld generators of the quantum affine superalgebra $U_q(\widehat{sl}(M|N))$.
In Section \ref{Section:3} we introduce bosons and 
give a bosonization of the quantum affine superalgebra $U_q(\widehat{sl}(M|N))$ for an arbitrary level $k \in {\bf C}$.
We realize the Wakimoto module as a submodule of this bosonization using the $\xi-\eta$ system.
In Section \ref{Section:4} we give screening operators that commute with $U_q(\widehat{sl}(M|N))$ for the level $k \neq -M+N$.
In Section \ref{Section:5} we prove the main results.
In Section \ref{Section:6} we give concluding remarks.
In Appendix \ref{appendix:A} we summarize normal ordering rules of bosonic operators.
In Appendix \ref{appendix:B} we recall a $q$-difference realization of $U_q(sl(M|N))$.
In Appendix \ref{appendix:C} we summarize useful formulae to take the limit $q\to 1$.

\section{$U_q(\widehat{sl}(M|N))$}
\label{Section:2}

In this Section we recall the definition of the quantum affine superalgebra $U_q(\widehat{sl}(M|N))$.

\subsection{Chevalley generator}

Throughout this paper, $q \in {\bf C}$ is assumed to be $0<|q|<1$.
For any integer $n$, define
\begin{eqnarray}
[n]_q=\frac{q^n-q^{-n}}{q-q^{-1}}.
\end{eqnarray}
We begin with the definition of
the quantum affine superalgebra $U_q(\widehat{sl}(M|N))$ for 
$M,N=1,2,3,\cdots$ in terms of
Chevalley generators.
The Cartan matrix of the affine Lie superalgebra $\widehat{sl}(M|N)$ is
\begin{eqnarray}
(A_{i,j})_{0\leq i,j \leq M+N-1}=
\left(\begin{array}{ccccccccccc}
0&-1&0&\cdots & & & & &\cdots &0&1\\
-1&2&-1&\cdots & & & & &\cdots &0&0\\
0&-1&2&\cdots & & & & &\cdots &\cdots &\cdots \\
\cdots&\cdots&\cdots&\cdots &-1&\cdots & & & & & \\
 & &\cdots&-1&2&-1&\cdots& & & & \\
 & & &\cdots &-1&0&1&\cdots& & & \\
 & & & & &\cdots 1&-2&1&\cdots & & \\
 & & & & &\cdots &1 &\cdots&\cdots&\cdots &\cdots\\
\cdots &\cdots &\cdots & & & & &\cdots&-2& 1&0\\
0&0&\cdots & & & & &\cdots&1&-2&1\\
1&0&\cdots & & & & &\cdots &0&1&-2 
\end{array}\right).
\end{eqnarray}
where the diagonal part is $(A_{i,i})_{0\leq i \leq M+N-1}=(0,\overbrace{2,2,\cdots,2}^{M-1},0,\overbrace{-2,-2,\cdots,-2}^{N-1})$.
Let $\bar{\alpha}_i$, $\bar{\Lambda}_i$ $(1\leq i \leq M+N-1)$ be the classical simple roots,
the classical fundamental weights, respectively.
Let $(\cdot|\cdot)$ be the symmetric bilinear form satisfying
$(\bar{\alpha}_i|\bar{\alpha}_j)=A_{i,j}$ and $(\bar{\Lambda}_i|\bar{\alpha}_j)=\delta_{i,j}$ for $1\leq i,j \leq M+N-1$.
Let us introduce the affine weight $\Lambda_0$ and the null root $\delta$ 
satisfying $(\Lambda_0|\Lambda_0)=(\delta|\delta)=0$, $(\Lambda_0|\delta)=1$,
and $(\Lambda_0|\bar{\alpha}_i)=(\Lambda_0|\bar{\Lambda}_i)=0$ for $1\leq i \leq M+N-1$.
The other affine weights and the affine roots are given by
$\Lambda_i=\bar{\Lambda}_i+\Lambda_0$, $\alpha_i=\bar{\alpha}_i$ for $1 \leq i \leq M+N-1$,
and $\alpha_0=\delta-\sum_{i=1}^{M+N-1}\alpha_i$.

The quantum affine superalgebra $U_q(\widehat{sl}(M|N))$ \cite{Yamane}
is the associative algebra over ${\bf C}$
with the Chevalley generators $e_i, f_i, h_i, d$ $(i=0,1,2,\cdots,M+N-1)$.
The ${\bf Z}_2$-grading of the Chevalley generators is given by
$p(e_0) \equiv p(f_0) \equiv p(e_{M}) \equiv p(f_{M})\equiv 1\pmod{2}$ and zero otherwise.
The defining relations of the Chevalley generators are given as follows.
\begin{eqnarray}
&&[h_i,h_j]=0,\\
&&[h_i,e_j]=A_{i,j}e_j,~~~[h_i,f_j]=-A_{i,j}f_j,\\
&&[e_i,f_j]=\delta_{i,j}\frac{q^{h_i}-q^{-h_i}}{q-q^{-1}},\\
&&[e_j,[e_j,e_i]_{q^{-1}}]_q=0,~~~[f_j,[f_j,f_i]_{q^{-1}}]_q=0~~~{\rm for}~|A_{i,j}|=1, i\neq 0, M,\\
&&[e_i,e_j]=0,~~~[f_i,f_j]=0~~~{\rm for}~|A_{i,j}|=0,
\\
&&[e_M,[e_{M+1},[e_M,e_{M-1}]_{q^{-1}}]_q]=0,~~~
[f_M,[f_{M+1},[f_M,f_{M-1}]_{q^{-1}}]_q]=0,\\
&&
[e_0,[e_1,[e_0,e_{M+N-1}]_{q}]_{q^{-1}}]=0,~~~
[f_0,[f_1,[f_0,f_{M+N-1}]_{q}]_{q^{-1}}]=0,
\end{eqnarray}
where we use the notation
\begin{eqnarray}
[X,Y]_a=XY-(-1)^{p(X)p(Y)}a YX,
\end{eqnarray}
for homogeneous elements $X,Y \in U_q(\widehat{sl}(M|N))$.
For simplicity we write $[X,Y]=[X,Y]_1$.
If $M=1$ or $N=1$, we have extra fifth order Serre relations.
As for the explicit forms of the extra Serre relations, we refer 
the reader to \cite{Koj3, Yamane}. 
Moreover, $U_q(\widehat{sl}(M|N))$ is a Hopf algebra over ${\bf C}$ with coproduct
\begin{eqnarray}
&&\Delta(h_i)=h_i \otimes 1+1 \otimes h_i,\\
&&\Delta(e_i)=e_i \otimes q^{h_i}+1 \otimes e_i,\\
&&\Delta(f_i)=f_i \otimes 1+q^{-h_i} \otimes f_i,
\end{eqnarray}
and antipode
\begin{eqnarray}
S(h_i)=-h_i,~~~
S(e_i)=-e_i q^{-h_i},~~~
S(f_i)=-q^{h_i}f_i.
\end{eqnarray}
The multiplication rule for the tensor product is ${\bf Z}_2$-graded and is defined for homogeneous elements 
$X_1, X_2, Y_1, Y_2 \in U_q(\widehat{sl}(M|N))$ by
$(X_1 \otimes Y_1) (X_2 \otimes Y_2)=(-1)^{p(Y_1)p(X_2)} (X_1 X_2 \otimes Y_1 Y_2)$,
which extends to inhomogeneous elements through linearity.
The coproduct is an algebra automorphism $\Delta(XY)=\Delta(X)\Delta(Y)$ 
and the antipode $S$ is a graded algebra anti-automorphism $S(XY)=(-1)^{p(X)p(Y)}S(Y)S(X)$.

\subsection{Drinfeld generator}

In \cite{Yamane} was given 
the second realization of the quantum affine superalgebra $U_q(\widehat{sl}(M|N))$ which is more convenient
for the concrete realization given in this article. 
We recall this realization that we call the Drinfeld realization \cite{Drinfeld}.
The quantum affine superalgebra $U_q(\widehat{sl}(M|N))$ is isomorphic to 
the associative algebra over ${\bf C}$ with the Drinfeld generators
$X_{m}^{\pm, i}~(i=1,2,\cdots,M+N-1, m \in {\bf Z})$,
$H_{n}^i~(i=1,2,\cdots,M+N-1, n \in {\bf Z}_{\neq 0})$,
$H^i~(i=1,2,\cdots,M+N-1)$, and $c$.
The ${\bf Z}_2$-grading of the Drinfeld generators is given by
$p(X_m^{\pm, M})\equiv 1 \pmod{2}$ for $m \in {\bf Z}$ and zero otherwise.
The defining relations of the Drinfeld generators are given as follows.
\begin{eqnarray}
&&c : {\rm central~element},\label{def:Drinfeld1}\\
&&[H^i,H_m^j]=0,\label{def:Drinfeld2}\\
&&[H_{m}^i,H_{n}^j]=\frac{[A_{i,j}m]_q[cm]_q}{m}\delta_{m+n,0},
\label{def:Drinfeld3}\\
&&[H^i,X^{\pm,j}(z)]=\pm A_{i,j}X^{\pm,j}(z),
\label{def:Drinfeld4}\\
&&[H_{m}^i, X^{\pm,j}(z)]=\pm \frac{[A_{i,j}m]_q}{m}q^{\mp \frac{c}{2}|m|} z^m X^{\pm,j}(z),
\label{def:Drinfeld5}
\\
&&(z_1-q^{\pm A_{i,j}}z_2)
X^{\pm,i}(z_1)X^{\pm,j}(z_2)
=
(q^{\pm A_{j,i}}z_1-z_2)
X^{\pm,j}(z_2)X^{\pm,i}(z_1),~~{\rm for}~|A_{i,j}|\neq 0,
\label{def:Drinfeld6}
\\
&&
[X^{\pm,i}(z_1), X^{\pm,j}(z_2)]=0,~~{\rm for}~|A_{i,j}|=0,
\label{def:Drinfeld7}
\\
&&[X^{+,i}(z_1), X^{-,j}(z_2)]
=\frac{\delta_{i,j}}{(q-q^{-1})z_1z_2}
\left(
\delta(q^{c}z_2/z_1)\Psi_+^i(q^{\frac{c}{2}}z_2)-
\delta(q^{-c}z_2/z_1)\Psi_-^i(q^{-\frac{c}{2}}z_2)\right), \label{def:Drinfeld8}
\\
&& 
[X^{\pm,i}(z_{1}),
[X^{\pm,i}(z_{2}), X^{\pm,j}(z)]_{q^{-1}}]_q+\left(z_1 \leftrightarrow z_2\right)=0,
~~~{\rm for}~|A_{i,j}|=1,~i\neq M,
\label{def:Drinfeld9}\\
&&
[X^{\pm,M}(z_1), [X^{\pm,M+1}(w_1), [X^{\pm,M}(z_2), X^{\pm, M-1}(w_2)]_{q^{-1}}]_q ]
+(z_1 \leftrightarrow z_2)=0,
\label{def:Drinfeld10}
\end{eqnarray}
where we set
$\delta(z)=\sum_{m \in {\bf Z}}z^m$.
Here we use 
the generating functions
\begin{eqnarray}
X^{\pm,j}(z)&=&
\sum_{m \in {\bf Z}}X_{m}^{\pm,j} z^{-m-1},\\
\Psi_\pm^i(q^{\pm \frac{c}{2}}z)&=&q^{\pm h_i}
\exp\left(
\pm (q-q^{-1})\sum_{m>0}H_{\pm m}^i z^{\mp m}
\right).
\end{eqnarray}
The Chevalley generators are obtained by
\begin{eqnarray}
h_i&=&H^i~~~(i=1,2,\cdots,M+N-1),\\
e_i&=&X_{0}^{+,i},~~~f_i=X_0^{-,i}~~~(i=1,2,\cdots,M+N-1),\\
h_0&=&c-(H^1+H^2+\cdots +H^{M+N-1}),\\
e_0&=&(-1)^N [X_0^{-,M+N-1},\cdots ,[X_0^{-,M+1}, [X_0^{-,M},\cdots,[X_0^{-,2},X_0^{-,1}]_{q^{-1}} \cdots ]_{q^{-1}}]_q \cdots ]_q,\\
f_0&=&q^{H^1+H^2+\cdots+H^{M+N-1}}\nonumber\\
&\times&
[\cdots[[\cdots[X_{-1}^{+,1},X_0^{+,2}]_q,\cdots, X_{0}^{+,M}]_q, X_0^{+, M+1}]_{q^{-1}},\cdots, X_0^{+,M+N-1}]_{q^{-1}}.
\end{eqnarray}

Let $V(\lambda)$ be the highest-weight module over $U_q(\widehat{sl}(M|N))$ generated by the highest weight vector $|\lambda \rangle \neq 0$ such that
\begin{eqnarray}
&&H_m^i |\lambda \rangle=X_m^{\pm,i}|\lambda \rangle=0~~~(m>0),\\
&&X_0^{+,i}|\lambda \rangle=0,\\
&&H^i|\lambda \rangle=l_i |\lambda \rangle,
\end{eqnarray}
where the classical part of the highest weight is $\bar{\lambda}=\sum_{i=1}^{M+N-1} l_i \bar{\Lambda}_i$.

\section{Bosonization of $U_q(\widehat{sl}(M|N))$}
\label{Section:3}

In this Section we give a bosonization of $U_q(\widehat{sl}(M|N))$ for an arbitrary level $k \in {\bf C}$.

\subsection{Boson}
In order to construct a bosonization of $U_q(\widehat{sl}(M|N))$,
we introduce bosons 
$a_m^i~(m \in {\bf Z}, 1\leq i \leq M+N-1)$,
$b_m^{i,j}~(m \in {\bf Z}, 1\leq i<j \leq M+N)$,
$c_m^{i,j}~(m \in {\bf Z}, \nu_i \nu_j=+1, 1\leq i<j \leq M+N)$,
and zero mode operators
$Q_a^i~(1\leq i \leq M+N-1)$,
$Q_b^{i,j}~(1\leq i<j \leq M+N)$,
$Q_c^{i,j}~(\nu_i \nu_j=+1, 1\leq i<j \leq M+N)$.
Their commutation relations are
\begin{eqnarray}
&&[a_m^i,a_n^j]=\frac{1}{m}[(k+g)m]_q[A_{i,j}m]_q\delta_{m+n,0},
~~~[a_0^i,Q_a^j]=(k+g)A_{i,j},\\
&&[b_m^{i,j},b_n^{i',j'}]=-\nu_i \nu_j \frac{1}{m} [m]_q^2 \delta_{i,i'}\delta_{j,j'}\delta_{m+n,0},
~~~[b_0^{i,j},Q_b^{i',j'}]=-\nu_i \nu_j \delta_{i,i'} \delta_{j,j'},\\
&&[c_m^{i,j},c_n^{i',j'}]=\nu_i \nu_j \frac{1}{m} [m]_q^2 \delta_{i,i'}\delta_{j,j'}\delta_{m+n,0},
~~~[c_0^{i,j},Q_c^{i',j'}]=\nu_i \nu_j \delta_{i,i'} \delta_{j,j'}~~~(\nu_i\nu_j=+1),\\
&&[Q_b^{i,j},Q_b^{i',j'}]=\pi \sqrt{-1}~~~(\nu_i \nu_j=\nu_{i'}\nu_{j'}=-1).
\end{eqnarray}
The remaining commutators vanish.
Here $g=M-N$
stands for the dual Coxeter number of $sl(M|N)$,
and
$\nu_i=1$ for $1\leq i \leq M$ and $\nu_i=-1$ for $M+1\leq i \leq M+N$.
For instance, we have a minus sigh by exchanging operators $e^{\pm Q_b^{i,j}}$,
\begin{eqnarray}
e^{\epsilon Q_b^{i,j}} e^{\epsilon' Q_b^{i',j'}}=- e^{\epsilon' Q_b^{i',j'}} e^{\epsilon Q_b^{i,j}}~~~
(\nu_i \nu_j=\nu_{i'}\nu_{j'}=-1, \epsilon,\epsilon'=\pm).
\end{eqnarray}
We define free boson fields $a_\pm^i(z), b_\pm^{i,j}(z), b^{i,j}(z), c^{i,j}(z)$ as follows.
\begin{eqnarray}
&&a_\pm^i(z)=\pm (q-q^{-1}) \sum_{m>0} a_{\pm m}^i z^{\mp m}\pm a_0^i {\rm log}q,
\\
&&b_\pm^{i,j}(z)=
\pm (q-q^{-1}) \sum_{m>0} b_{\pm m}^{i,j} z^{\mp m}\pm b_0^{i,j} {\rm log}q,
\\
&&b^{i,j}(z)=-\sum_{m \neq 0}\frac{b_m^{i,j}}{[m]_q}z^{-m}+Q_b^{i,j}+b_0^{i,j}{\rm log}z,\\
&&c^{i,j}(z)=-\sum_{m \neq 0}\frac{c_m^{i,j}}{[m]_q}z^{-m}+Q_c^{i,j}+c_0^{i,j}{\rm log}z.
\end{eqnarray}
We define free boson fields $(\Delta^{\varepsilon}_{L} b_\pm^{i,j})(z),
(\Delta^{\varepsilon}_{R} b_\pm^{i,j})(z)$ $(\varepsilon=\pm,0)$ as follows.
\begin{eqnarray}
&&(\Delta^{\varepsilon}_L b_\pm^{i.j})(z)=
\left\{\begin{array}{cc}
b_\pm^{i+1,j}(q^{\varepsilon}z)-b_\pm^{i,j}(z)& (\varepsilon=\pm),\\
b_\pm^{i+1,j}(z)+b_\pm^{i,j}(z)& (\varepsilon=0),
\end{array}\right.
\\
&&(\Delta^{\varepsilon}_R b_\pm^{i.j})(z)=\left\{\begin{array}{cc}
b_\pm^{i,j+1}(q^{\varepsilon}z)-b_\pm^{i,j}(z)& (\varepsilon=\pm),\\
b_\pm^{i,j+1}(z)+b_\pm^{i,j}(z)& (\varepsilon=0).
\end{array}\right.
\end{eqnarray}
We define further
free boson fields with parameters $L_1,L_2,\cdots,L_r, M_1,M_2,\cdots, M_r, \alpha$ as follows.
\begin{eqnarray}
&&
\left(\frac{L_1}{M_1}\frac{L_2}{M_2}\cdots \frac{L_r}{M_r}~a^i\right)(z;\alpha)\nonumber
\\
&=&-\sum_{m \neq 0}
\frac{[L_1 m]_q [L_2 m]_q \cdots [L_r m]_q}{[M_1 m]_q [M_2m]_q \cdots [M_r m]_q}\frac{a_m^i}{[m]_q}q^{-\alpha|m|}z^{-m}
+\frac{L_1 L_2 \cdots L_r}{M_1 M_2 \cdots M_r}(Q_a^i+a_0^i {\rm log}z).
\end{eqnarray} 
Normal ordering rules are defined as follows.
\begin{eqnarray}
&&:a_m^i a_n^j:=:a_n^j a_m^i:=\left\{
\begin{array}{cc}
a_m^i a_n^j& (m<0),\\
a_n^j a_m^i& (m>0),
\end{array}
\right.
\\
&&:b_m^{i,j} b_n^{i',j'}:=:b_n^{i',j'} b_m^{i,j}:=\left\{
\begin{array}{cc}
b_m^{i,j} b_n^{i',j'}& (m<0),\\
b_n^{i',j'} b_m^{i,j}& (m>0),
\end{array}
\right.
\\
&&:Q_b^{i,j} Q_b^{i',j'}:=:Q_b^{i',j'} Q_b^{i,j}:=Q_b^{i,j} Q_b^{i',j'}~~~(i>i'~~{\rm or}~~i=i',j>j').
\end{eqnarray}
Normal ordering rules of $c_m^{i,j}$ and $Q_c^{i,j}$ are defined in the same way.
For instance we have
\begin{eqnarray}
&&:\exp(a^i(z)):=\exp\left(\sum_{m>0}\frac{a_{-m}^i}{[m]_q}z^{m}\right)
\exp\left(-\sum_{m>0}\frac{a_m^i}{[m]_q}z^{-m}\right)e^{Q_a^i}z^{a_0^i},\\
&&:e^{Q_b^{i,j}}e^{Q_b^{i',j'}}:=:e^{Q_b^{i',j'}}e^{Q_b^{i,j}}:=
e^{Q_b^{i,j}}e^{Q_b^{i',j'}}=\nu_i \nu_j e^{Q_b^{i',j'}}e^{Q_b^{i,j}}~~~
(i>i'~{\rm or}~i=i',j>j').
\end{eqnarray}

\subsection{Bosonization}

In this Section we fix a complex number $k \in {\bf C}$.\\
$\bullet$~We define bosonic operators $H^i(z)~(1\leq i \leq M+N-1)$ as follows.
\begin{eqnarray}
H^i(z)&=&\frac{1}{(q-q^{-1})z}\left\{
\left[a_+^i(q^{\frac{g}{2}}z)
+\sum_{l=1}^i (\Delta_{R}^- b_+^{l,i})(q^{\frac{k}{2}+l}z)\right.\right.\nonumber\\
&&\left.\left.-\sum_{l=i+1}^M (\Delta_{L}^{-} b_+^{i,l})(q^{\frac{k}{2}+l}z)
-\sum_{l=M+1}^{M+N}
(\Delta_L^-  b_+^{i,l})(q^{\frac{k}{2}+2M+1-l}z)
\right]\right.
\nonumber\\
&&\left.-\left[a_-^i(q^{-\frac{g}{2}}z)
+\sum_{l=1}^i (\Delta_{R}^+ b_-^{l,i})(q^{-\frac{k}{2}-l}z)\right.\right.\nonumber
\\
&&\left.\left.
-\sum_{l=i+1}^M (\Delta_{L}^{+} b_-^{i,l})(q^{-\frac{k}{2}-l}z)
-\sum_{l=M+1}^{M+N}
(\Delta_L^+  b_-^{i,l})(q^{-\frac{k}{2}-2M-1+l}z)
\right]\right\}
~~~(1\leq i \leq M-1),
\label{def:H1}
\\
H^M(z)
&=&\frac{1}{(q-q^{-1})z}\left\{
\left[a_+^M(q^{\frac{g}{2}}z)-
\sum_{l=1}^{M-1} (\Delta_{R}^0 b_+^{l,M})(q^{\frac{k}{2}+l}z)\right.\right.
\left.+\sum_{l=M+2}^{M+N} (\Delta_L^0 b_+^{M,l})(q^{\frac{k}{2}+2M+1-l}z)
\right]\nonumber\\
&&-
\left[a_-^M(q^{-\frac{g}{2}}z)-
\sum_{l=1}^{M-1} (\Delta_R^0 b_-^{l,M})(q^{-\frac{k}{2}-l}z)\right.
\left.\left.+\sum_{l=M+2}^{M+N} (\Delta_L^0 b_-^{M,l})(q^{-\frac{k}{2}-2M-1+l}z)
\right]\right\},
\label{def:H2}
\\
H^i(z)
&=&\frac{1}{(q-q^{-1})z}\left\{
\left[a_+^i(q^{\frac{g}{2}}z)
-\sum_{l=1}^M (\Delta_R^+ b_+^{l,i})(q^{\frac{k}{2}+l-1}z)\right.\right.\nonumber\\
&&\left.-\sum_{l=M+1}^i (\Delta_R^+ b_+^{l,i})(q^{\frac{k}{2}+2M-l}z)
+\sum_{l=i+1}^{M+N}
(\Delta_L^+ b_+^{i,l})(q^{\frac{k}{2}+2M-l}z)
\right]\nonumber\\
&&-\left[a_-^i(q^{-\frac{g}{2}}z)
-\sum_{l=1}^M (\Delta_R^- b_-^{l,i})(q^{-\frac{k}{2}-l+1}z)\right.
-\sum_{l=M+1}^i (\Delta_R^- b_-^{l,i})(q^{-\frac{k}{2}-2M+l}z)\nonumber
\\
&&\left.\left.+\sum_{l=i+1}^{M+N}
(\Delta_L^- b_-^{i,l})(q^{-\frac{k}{2}-2M+l}z)
\right]\right\}~~~(M+1\leq i \leq M+N-1).
\label{def:H3}
\end{eqnarray}
We define bosonic operators $H_m^i~(m \in {\bf Z}, 1\leq i \leq M+N-1)$ as follows.
\begin{eqnarray}
H^i(z)=\sum_{m \in {\bf Z}}H_m^i z^{-m-1}~~~(1\leq i \leq M+N-1).\label{def:H4}
\end{eqnarray}
$\bullet$~We define bosonic operators $\Psi_\pm^i(z)~(1\leq i \leq M+N-1)$ as follows.
\begin{eqnarray}
\Psi_\pm^i(q^{\pm \frac{k}{2}}z)
&=&
:e^{a_\pm^i(q^{\pm \frac{g}{2}}z)+\sum_{l=1}^i (\Delta_R^{\mp}b_\pm^{l,i})(q^{\pm(\frac{k}{2}+l)}z)
-\sum_{l=i+1}^M (\Delta_L^{\mp}b_\pm^{i,l})(q^{\pm(\frac{k}{2}+l)}z)}\nonumber\\
&\times& 
e^{-\sum_{l=M+1}^{M+N}(\Delta_L^{\mp} b_\pm^{i,l})(q^{\pm(\frac{k}{2}+2M+1-l)}z)}:~~~(1\leq i \leq M-1),
\label{def:psi1}
\\
\Psi_\pm^M(q^{\pm \frac{k}{2}}z)
&=&
:e^{a_\pm^M(q^{\pm \frac{g}{2}}z)-\sum_{l=1}^{M-1} (\Delta_R^0 b_\pm^{l,M})(q^{\pm(\frac{k}{2}+l)}z)+\sum_{l=M+2}^{M+N} (\Delta_L^0 b_\pm^{M,l})(q^{\pm(\frac{k}{2}+2M+1-l)}z)}:,
\label{def:psi2}
\\
\Psi_\pm^i(q^{\pm \frac{k}{2}}z)
&=&:e^{a_\pm^i(q^{\pm \frac{g}{2}}z)
-\sum_{l=1}^M (\Delta_R^\pm b_\pm^{l,i})(q^{\pm (\frac{k}{2}+l-1)}z)-\sum_{l=M+1}^i (\Delta_R^\pm b_\pm^{l,i})(q^{\pm(\frac{k}{2}+2M-l)}z)}\nonumber\\
&\times&
e^{\sum_{l=i+1}^{M+N}
(\Delta_L^\pm  b_\pm^{i,l})(q^{\pm(\frac{k}{2}+2M-l)}z)}:~~~(M+1\leq i \leq M+N-1).
\label{def:psi3}
\end{eqnarray}
$\bullet$~We define bosonic operators $X^{+,i}(z)~(1\leq i \leq M+N-1)$ as follows.
\begin{eqnarray}
X^{+,i}(z)&=&\sum_{j=1}^i\frac{c_{i,j}}{(q-q^{-1})z}(E_{i,j}^+(z)-E_{i,j}^-(z))~~~(1\leq i \leq M-1),
\label{def:X^+1}\\
X^{+,M}(z)&=&\sum_{j=1}^M c_{M,j} E_{M,j}(z),
\label{def:X^+2}
\\
X^{+,i}(z)&=&
\sum_{j=1}^M c_{i,j} E_{i,j}(z)+\sum_{j=M+1}^i \frac{c_{i,j}}{(q-q^{-1})z}(E_{i,j}^+(z)-E_{i,j}^-(z))~~~~~
(M+1\leq i \leq M+N-1).\nonumber\\
\label{def:X^+3}
\end{eqnarray}
Here $c_{i,j} \in {\bf C}_{\neq 0}~(1\leq j \leq i \leq M+N-1)$ and we set $E_{i,j}^\pm(z)$ and $E_{i,j}(z)$ as follows.
\\
For $1\leq i \leq M-1$ and $1\leq j \leq i-1$ we set
\begin{eqnarray}
E_{i,j}^\pm (z)&=&:e^{(b+c)^{j,i}(q^{j-1}z)+b_\pm^{j,i+1}(q^{j-1}z)-(b+c)^{j,i+1}(q^{j-1\pm 1}z)
+\sum_{l=1}^{j-1}(\Delta_R^- b_+^{l,i})(q^lz)}:,
\\
E_{i,i}^\pm (z)
&=&:e^{b_\pm^{i,i+1}(q^{i-1}z)-(b+c)^{i,i+1}(q^{i-1\pm 1}z)+\sum_{l=1}^{i-1}
(\Delta_R^- b_+^{l,i})(q^lz)}:.
\end{eqnarray}
For $1\leq j \leq M-1$ we set
\begin{eqnarray}
E_{M,j}(z)
&=&:e^{(b+c)^{j,M}(q^{j-1}z)+b^{j,M+1}(q^{j-1}z)-\sum_{l=1}^{j-1}(\Delta_R^0 b_+^{l,M})(q^lz)}:,
\\
E_{M,M}(z)
&=&:e^{b^{M,M+1}(q^{M-1}z)-\sum_{l=1}^{M-1}(\Delta_R^0 b_+^{l,M})(q^lz)}:.
\end{eqnarray}
For $M+1 \leq i \leq M+N-1$ we set
\begin{eqnarray}
E_{i,i}^\pm (z)
&=&:e^{-b_\pm ^{i,i+1}(q^{2M+1-i}z)-(b+c)^{i,i+1}(q^{2M+1 \mp 1-i}z)}\nonumber\\
&\times&
e^{-\sum_{l=1}^M(\Delta_R^+ b_+^{l,i})(q^{l-1}z)
-\sum_{l=M+1}^{i-1}(\Delta_R^+ b_+^{l,i})(q^{2M-l}z)}:.
\end{eqnarray}
For $M+1 \leq i \leq M+N-1$ and $1\leq j \leq M$ we set
\begin{eqnarray}
E_{i,j}(z)
&=&:e^{b_+^{j,i}(q^{j-1}z)-b^{j,i}(q^jz)+b^{j,i+1}(q^{j-1}z)
-\sum_{l=1}^{j-1}(\Delta_R^+ b_+^{l,i})(q^{l-1}z)}:.
\end{eqnarray}
For $M+1 \leq i \leq M+N-1$ and $M+1\leq j \leq i-1$ we set
\begin{eqnarray}
E_{i,j}^\pm (z)&=&:e^{(b+c)^{j,i}(q^{2M+1-j}z)-b_\pm^{j,i+1}(q^{2M+1-j}z)-(b+c)^{j,i+1}(q^{2M+1 \mp 1-j}z)}\nonumber\\
&\times&
e^{-\sum_{l=1}^M(\Delta_R^+ b_+^{l,i})(q^{l-1}z)-\sum_{l=M+1}^{j-1}(\Delta_R^+ b_+^{l,i})(q^{2M-l}z)}:.
\end{eqnarray}
$\bullet$~We define bosonic operators $X^{-,i}(z)~(1\leq i \leq M+N-1)$ as follows.
\begin{eqnarray}
X^{-,i}(z)&=&
\sum_{j=1}^{i-1}
\frac{d_{i,j}^1}{(q-q^{-1})z}(F_{i,j}^{1,-}(z)-F_{i,j}^{1,+}(z))
+\frac{d_{i,i}^2}{(q-q^{-1})z}(F_{i,i}^{2,-}(z)-F_{i,i}^{2,+}(z))\nonumber\\
&+&\sum_{j=i+2}^M \frac{d_{i,j}^3}{(q-q^{-1})z}(F_{i,j}^{3,-}(z)-F_{i,j}^{3,+}(z))+\sum_{j=M+1}^{M+N} d_{i,j}^3 F_{i,j}^3(z)~~~~~(1\leq i \leq M-1),
\label{def:X^-1}
\\
X^{-,M}(z)
&=&\sum_{j=1}^{M-1}\frac{d_{M,j}^1}{(q-q^{-1})z}(F_{M,j}^{1,-}(z)-F_{M,j}^{1,+}(z))+\frac{d_{M,M}^2}{(q-q^{-1})z}(F_{M,M}^{2,-}(z)-F_{M,M}^{2,+}(z))\nonumber\\
&+&\sum_{j=M+2}^{M+N}\frac{d_{M,j}^3}{(q-q^{-1})z}(F_{M,j}^{3,-}(z)-F_{M,j}^{3,+}(z)),
\label{def:X^-2}
\\
X^{-,i}(z)
&=&\sum_{j=1}^M d_{i,j}^1 F_{i,j}^{1}(z)
+\sum_{j=M+1}^{i-1}\frac{d_{i,j}^1}{(q-q^{-1})z}(F_{i,j}^{1,-}(z)-F_{i,j}^{1,+}(z))+
\frac{d_{i,i}^2}{(q-q^{-1})z}(F_{i,i}^{2,-}(z)-F_{i,i}^{2,+}(z))\nonumber
\\
&&+\sum_{j=i+2}^{M+N}\frac{d_{i,j}^3}{(q-q^{-1})z}(F_{i,j}^{3,-}(z)-F_{i,j}^{3,+}(z))~~~~~(M+1 \leq i \leq M+N-1).
\label{def:X^-3}
\end{eqnarray}
Here we set $F_{i,j}^{1,\pm}(z), F_{i,j}^1(z)$, $F_{i,i}^{2,\pm}(z)$, $F_{i,j}^{3,\pm}(z), F_{i,j}^3(z)$ as follows.
\\
For $1\leq i \leq M-1$ and $1\leq j \leq i-1$ we set
\begin{eqnarray}
F_{i,j}^{1,\pm}(z)
&=&:e^{a_-^i(q^{-\frac{k+g}{2}}z)+(b+c)^{j,i+1}(q^{-k-j}z)-b_\pm^{j,i}(q^{-k-j}z)-(b+c)^{j,i}(q^{-k-j\mp 1}z)}\nonumber\\
&\times& e^{\sum_{l=j+1}^i (\Delta_R^+ b_-^{l,i})(q^{-k-l}z)-\sum_{l=i+1}^M (\Delta_L^+ b_-^{i,l})(q^{-k-l}z)
-\sum_{l=M+1}^{M+N}(\Delta_L^+ b_-^{i,l})(q^{-k-2M-1+l}z)}:.
\end{eqnarray}
For $1\leq j \leq M-1$ we set
\begin{eqnarray}
F_{M,j}^{1,\pm}(z)&=&:e^{a_-^M(q^{-\frac{k+g}{2}}z)-b_\pm^{j,M}(q^{-k-j}z)-(b+c)^{j,M}(q^{-k-j \mp 1}z)-b_-^{j,M+1}(q^{-k-j}z)-b^{j,M+1}(q^{-k-j+1}z)}
\nonumber\\
&\times& e^{-\sum_{l=j+1}^{M-1}(\Delta_R^0 b_-^{l,M})(q^{-k-l}z)
+\sum_{l=M+2}^{M+N}(\Delta_L^0 b_-^{M,l})(q^{-k-2M-1+l}z)}:.
\end{eqnarray}
For $M+1 \leq i \leq M+N-1$ and $1\leq j \leq M$ we set
\begin{eqnarray}
F_{i,j}^1(z)&=&
:e^{a_-^i(q^{-\frac{k+g}{2}}z)-b_-^{j,i+1}(q^{-k-j}z)-b^{j,i+1}(q^{-k-j+1}z)+b^{j,i}(q^{-k-j}z)
-\sum_{l=j+1}^M(\Delta_R^- b_-^{l,i})(q^{-k-l+1}z)}
\nonumber\\
&\times&
e^{-\sum_{l=M+1}^i (\Delta_R^- b_-^{l,i})(q^{-k-2M+l}z)
+\sum_{l=i+1}^{M+N}(\Delta_L^- b_-^{i,l})(q^{-k-2M+l}z)}:.
\end{eqnarray}
For $M+1 \leq i \leq M+N-1$ and $M+1\leq j \leq i-1$ we set
\begin{eqnarray}
F_{i,j}^{1,\pm}(z)&=&:e^{a_-^i(q^{-\frac{k+g}{2}}z)+(b+c)^{j,i+1}(q^{-k-2M+j}z)
+b_\pm^{j,i}(q^{-k-2M+j}z)-(b+c)^{j,i}(q^{-k-2M\pm 1+j}z)}\nonumber
\\
&\times& e^{-\sum_{l=j+1}^i(\Delta_R^- b_-^{l,i})(q^{-k-2M+l}z)
+\sum_{l=i+1}^{M+N}(\Delta_L^- b_-^{i,l})(q^{-k-2M+l}z)}:.
\end{eqnarray}
For $1 \leq i \leq M-1$ we set
\begin{eqnarray}
F_{i,i}^{2,\pm}(z)&=&:e^{a_\pm^i(q^{\pm \frac{k+g}{2}}z)
+b_\pm^{i,i+1}(q^{\pm(k+i+1)}z)+(b+c)^{i,i+1}(q^{\pm(k+i)}z)}
\nonumber\\
&\times&
e^{-\sum_{l=i+2}^M(\Delta_L^\mp b_\pm^{i,l})(q^{\pm (k+l)}z)
-\sum_{l=M+1}^{M+N}(\Delta_L^\mp b_\pm^{i,l})(q^{\pm (k+2M+1-l)}z)}:,
\\
F_{M,M}^{2,\pm}(z)
&=&:e^{a_\pm^M(q^{\pm \frac{k+g}{2}}z)-b^{M,M+1}(q^{\pm(k+M-1)}z)
+\sum_{l=M+2}^{M+N}(\Delta_L^0 b_\pm^{M,l})(q^{\pm(k+2M+1-l)}z)}:.
\end{eqnarray}
For $M+1 \leq i \leq M+N-1$ we set
\begin{eqnarray}
F_{i,i}^{2,\pm}(z)&=&:e^{a_\pm^i(q^{\pm \frac{k+g}{2}}z)
-b_\pm^{i,i+1}(q^{\pm(k+2M-1-i)}z)+(b+c)^{i,i+1}(q^{\pm (k+2M-i)}z)}\nonumber\\
&\times& 
e^{\sum_{l=i+2}^{M+N}(\Delta_L^\pm b_\pm^{i,l})(q^{\pm(k+2M-l)}z)}:.
\end{eqnarray}
For $1\leq i \leq M-2$ and $i+2 \leq j \leq M$ we set
\begin{eqnarray}
F_{i,j}^{3,\pm}(z)&=&:e^{a_+^i(q^{\frac{k+g}{2}}z)+
(b+c)^{i,j}(q^{k+j-1}z)
+b_\pm^{i+1,j}(q^{k+j-1}z)-(b+c)^{i+1,j}(q^{k-1\pm 1+j}z)}\nonumber\\
&\times&e^{-\sum_{l=j}^M(\Delta_L^- b_+^{i,l})(q^{k+l}z)
-\sum_{l=M+1}^{M+N}(\Delta_L^- b_+^{i,l})(q^{k+2M+1-l}z)}:.
\end{eqnarray}
For $1\leq i \leq M-1$ and $M+1 \leq j \leq M+N$ we set
\begin{eqnarray}
F_{i,j}^{3}(z)&=&:e^{a_+^i(q^{\frac{k+g}{2}}z)-b^{i,j}(q^{k+2M-j}z)-b_+^{i+1,j}(q^{k+2M-j}z)+b^{i+1,j}(q^{k+2M+1-j}z)}\nonumber\\
&\times&
e^{-\sum_{l=j+1}^{M+N}(\Delta_L^- b_+^{i,l})(q^{k+2M+1-l}z)}:.
\end{eqnarray}
For $M+2 \leq j \leq M+N$ we set
\begin{eqnarray}
F_{M,j}^{3,\pm}(z)&=&
:e^{a_+^M(q^{\frac{k+g}{2}}z)-b^{M,j}(q^{k+2M-j}z)-b_\pm^{M+1,j}(q^{k+2M+1-j}z)
-(b+c)^{M+1,j}(q^{k+2M+1 \mp 1-j}z)}\nonumber\\
&\times& e^{b_+^{M+1,j}(q^{k+2M+1-j}z)+\sum_{l=j+1}^{M+N}(\Delta_L^0 b_+^{M,l})(q^{k+2M+1-l}z)}:.
\end{eqnarray}
For $M+1\leq i \leq M+N-1$ and $i+2 \leq j \leq M+N$ we set
\begin{eqnarray}
F_{i,j}^{3,\pm}(z)&=&:e^{a_+^i(q^{\frac{k+g}{2}}z)+(b+c)^{i,j}(q^{k+2M+1-j}z)
-b_\pm^{i+1,j}(q^{k+2M+1-j}z)-(b+c)^{i+1,j}(q^{k+2M+1 \mp1-j}z)}\nonumber\\
&\times& e^{\sum_{l=j+1}^{M+N}(\Delta_L^+ b_+^{i,l})(q^{k+2M-l}z)}:.
\end{eqnarray}
Here we set $d_{i,j}^1, d_{i,i}^2, d_{i,j}^3 \in {\bf C}$ as follows.
\begin{eqnarray}
d_{i,j}^1&=&\nu_{i+1}\frac{1}{c_{i,j}}\times
\left\{
\begin{array}{cc}
1& (1 \leq i \leq M-1, 1 \leq j \leq i-1),
\\
q^{j-1} & (i=M, 1\leq j \leq M-1),
\\
q^{-k-1} & (M+1 \leq i \leq M+N-1, 1 \leq j \leq M),
\\
1 & (M+1 \leq i \leq M+N-1, M+1 \leq j \leq i-1),
\end{array}
\right.\label{def:d1}
\\
d_{i,i}^2&=&
\nu_{i+1}
\frac{1}{c_{i,i}}
\times
\left\{\begin{array}{cc}
1 & (1\leq i \neq M \leq M+N-1),\\
q^{M-1} & (i=M),
\end{array}
\right.\label{def:d2}
\\
d_{i,j}^3
&=&\nu_{i+1}
\frac{1}{c_{i,i}}\prod_{l=1}^{j-i-1}\frac{c_{i+l,i+1}}{c_{i+l,i}}
\times \left\{
\begin{array}{cc}
1& (1 \leq i \leq M-1, i+2 \leq j \leq M),\\
q^{k+3M+1-2j}& (1\leq i \leq M-1, M+1 \leq j \leq M+N),\\
q^{(M-1)(j-M)}& (i=M, M+2 \leq j \leq M+N),\\
1 & (M+1 \leq i \leq M+N-1, i+2 \leq j \leq M+N).
\end{array}\right.\label{def:d3}
\end{eqnarray}

The following is {\bf the
first main result} of this article.

\begin{thm}~~~The bosonic operators
$H^i=H_0^i,~H_m^i~(m \in {\bf Z}_{\neq 0}, 1\leq i \leq M+N-1)$ defined in (\ref{def:H1})-(\ref{def:H4}),
$\Psi_\pm^i(z)~(1\leq i \leq M+N-1)$ defined in (\ref{def:psi1})-(\ref{def:psi3}), and
$X^{\pm,i}(z)~(1\leq i \leq M+N-1)$ defined in (\ref{def:X^+1})-(\ref{def:X^+3}) 
and (\ref{def:X^-1})-(\ref{def:X^-3}) satisfy the defining relations of the Drinfeld realization 
(\ref{def:Drinfeld1})-(\ref{def:Drinfeld10}) with the central element $c=k \in {\bf C}$.
Here
the coefficients $d_{i,j}^1, d_{i,i}^2$, and $d_{i,j}^3$ 
are given
in (\ref{def:d1})-(\ref{def:d3}).
\label{thm:1}
\end{thm}
This bosonization reproduces those of $U_q(\widehat{sl}(M|1))$ 
upon the specialization $N=1$ \cite{Koj1}.

We introduce the boson Fock space ${F}(p_a,p_b,p_c)$ as follows.
The vacuum state $|0\rangle \neq 0$ is defined by
\begin{eqnarray}
a_m^i|0\rangle=
b_m^{i,j}|0\rangle=
c_m^{i,j}|0\rangle=0~~~(m \geq 0).
\end{eqnarray}
Let $|p_a,p_b,p_c\rangle$ be
\begin{eqnarray}
|p_a,p_b,p_c\rangle&=&\exp\left(
\sum_{i,j=1}^{M+N-1}\frac{(A^{-1})_{i,j}}{k+g}p_a^i Q_a^i-\sum_{1\leq i<j \leq M+N-1}\nu_i \nu_j p_b^{i,j}Q_b^{i,j}+\sum_{1\leq i<j \leq M+N-1
\atop{\nu_i \nu_j=+1}}p_c^{i,j}Q_c^{i,j}\right)|0\rangle,\nonumber\\
\end{eqnarray}
then $|p_a, p_b, p_c\rangle$
is the highest weight state of the boson Fock space ${F}(p_a,p_b,p_c)$, i.e.,
\begin{eqnarray}
&&
a_m^i|p_a,p_b,p_c \rangle=
b_m^{i,j}|p_a,p_b,p_c \rangle=
c_m^{i,j}|p_a,p_b,p_c \rangle=0~~~(m>0),
\\
&&
a_0^i|p_a,p_b,p_c \rangle=p_a^i |p_a,p_b,p_c\rangle,~~
b_0^{i,j}|p_a,p_b,p_c \rangle=p_b^{i,j}|p_a,p_b,p_c\rangle,\\
&&
c_0^{i,j}|p_a,p_b,p_c \rangle=p_c^{i,j}|p_a,p_b,p_c\rangle~~~(\nu_i \nu_j=+1).
\end{eqnarray}
The boson Fock space $F(p_a,p_b,p_c)$ is generated by the bosons $a_m^i, b_m^{i,j}, c_m^{i,j}$ on the highest weight state
$|p_a,p_b,p_c\rangle$.
We set the space $F(p_a)$ by
\begin{eqnarray}
F(p_a)=\bigoplus_{p_b^{i,j}=-p_c^{i,j} \in {\bf Z}~(\nu_i\nu_j=+)
\atop{p_b^{i,j}}\in {\bf Z}~(\nu_i\nu_j=-)}F(p_a,p_b,p_c).
\end{eqnarray}
Here we impose the restriction
$p_b^{i,j}=-p_c^{i,j}~(\nu_i \nu_j=+)$, because $X_m^{\pm,i}$ change $Q_b^{i,j}+Q_c^{i,j}$.
$F(p_a)$ is $U_q(\widehat{sl}(M|N))$-module.
Let $|\lambda \rangle=|p_a,0,0\rangle$ where $p_a^i=l_i~(1\leq i \leq M+N-1)$.

\begin{prop}~~~~The Drinfeld generators $H^i, H_m^i, X_m^{\pm, i}$ act on $|\lambda\rangle$ as follows.
\begin{eqnarray}
&&H_m^i |\lambda \rangle=X_m^{\pm,i}|\lambda \rangle=0~~~(m>0),\\
&&X_0^{+,i}|\lambda \rangle=0,\\
&&H^i|\lambda \rangle=l_i |\lambda \rangle.
\end{eqnarray}
\end{prop}
This property is just the highest weight state condition of the highest weight module $V(\lambda)$.

\begin{cor}~~~
We have the level-$k$ highest weight module $V(\lambda)$ of $U_q(\widehat{sl}(M|N))$ :
\begin{eqnarray}
V(\lambda) \subset F(p_a).
\end{eqnarray}
Here the classical part of the highest weight is $\bar{\lambda}=\sum_{i=1}^{M+N-1} l_i \bar{\Lambda}_i$.
\end{cor}

The module $F(p_a)$ is not irreducible.
In \cite{Kon} the irreducible highest weight module of $U_q(\widehat{sl}(2))$ was constructed by two steps
from the similar module as $F(p_a)$ :
the first step is construction of Wakimoto module using $\xi-\eta$ system,
and the second step is resolution by Felder complex using screening operators $Q_i$ \cite{Fel}.
Construction of Felder complex is an open problem even for non-superalgebra $U_q(\widehat{sl}(3))$.
In this paper we propose Wakimoto module of $U_q(\widehat{sl}(M|N))$ using $\xi-\eta$ system.
We would like to report on Felder-complex of $U_q(\widehat{sl}(M|N))$ in future publication.

We set bosonic operators $\xi_m^{i,j}, \eta_m^{i,j}~(\nu_i \nu_j=+1, 1\leq i<j \leq M+N-1)$ as follows.
\begin{eqnarray}
\eta^{i,j}(z)=\sum_{m \in {\bf Z}}\eta_m^{i,j} z^{-m-1}=:e^{c^{i,j}(z)}:,~~~
\xi^{i,j}(z)=\sum_{m \in {\bf Z}}\xi_m^{i,j}z^{-m}=:e^{-c^{i,j}(z)}:.
\end{eqnarray}
Fourier components $\eta_m^{i,j}=\oint \frac{dz}{2\pi \sqrt{-1}}z^m \eta^{i,j}(z),
\xi_m^{i,j}=\oint \frac{dz}{2\pi \sqrt{-1}}z^{m-1}\xi^{i,j}(z)$ are well-defined on the module $F(p_a)$.
The ${\bf Z}_2$-grading is given by
$p(\xi_m^{i,j})=p(\eta_m^{i,j})=+1$.
They satisfy anti-commutation relations.
\begin{eqnarray}
[\eta_m^{i,j},\xi_m^{i,j}]=\delta_{m+n,0},~~
[\eta_m^{i,j},\eta_n^{i,j}]=[\xi_m^{i,j},\xi_n^{i,j}]=0.
\end{eqnarray}
They commute with each other.
\begin{eqnarray}
[\eta_m^{i,j},\xi_m^{i',j'}]=
[\eta_m^{i,j},\eta_n^{i',j'}]=[\xi_m^{i,j},\xi_n^{i',j'}]=0~~~((i,j)\neq (i',j')).
\end{eqnarray}
The operators $\eta_0^{i,j}, \xi_0^{i,j}$ satisfy
\begin{eqnarray}
&&{\rm Im}(\eta_0^{i,j})={\rm Ker}(\eta_0^{i,j}),~~~
{\rm Im}(\xi_0^{i,j})={\rm Ker}(\xi_0^{i,j}),\\
&&\eta_0^{i,j} \xi_0^{i,j}+\xi_0^{i,j} \eta_0^{i,j}=1,\\
&&
(\eta_0^{i,j} \xi_0^{i,j})^2=\eta_0^{i,j}\xi_0^{i,j},~~~
(\xi_0^{i,j} \eta_0^{i,j})^2=\xi_0^{i,j}\eta_0^{i,j},\\
&&(\xi_0^{i,j}\eta_0^{i,j})(\eta_0^{i,j}\xi_0^{i,j})=
(\eta_0^{i,j}\xi_0^{i,j})(\xi_0^{i,j} \eta_0^{i,j})=0.
\end{eqnarray}
Hence we have direct sum decomposition.
\begin{eqnarray}
F(p_a)=\eta_0^{i,j}\xi_0^{i,j} F(p_a) \oplus \xi_0^{i,j} \eta_0^{i,j} F(p_a),
\end{eqnarray}
where
${\rm Ker}(\eta_0^{i,j})=\eta_0^{i,j} \xi_0^{i,j} F(p_a)$,
${\rm Coker}(\eta_0^{i,j})=\xi_0^{i,j}\eta_0^{i,j} F(p_a)$.
We set
\begin{eqnarray}
\eta_0=\prod_{1\leq i<j \leq M+N-1\atop{\nu_i\nu_j=+1}}\eta_0^{i,j},~~~
\xi_0=\prod_{1\leq i<j \leq M+N-1\atop{\nu_i\nu_j=+1}}\xi_0^{i,j}.
\end{eqnarray}
We introduce the subspace ${\cal F}(p_a)$ that gives a generalization of
the articles \cite{ZG, Koj2} by
\begin{eqnarray}
{\cal F}(p_a)=\eta_0\xi_0 F(p_a).
\end{eqnarray}
The operators $\eta_0^{i,j}, \xi_0^{i,j}$ commute with $X^{\pm,i}(z), \Psi_\pm^i(z)$ up to sign $\pm 1$.

\begin{prop}~~~${\cal F}(p_a)$ is $U_q(\widehat{sl}(M|N))$-module.
\end{prop}
We call ${\cal F}(p_a)$ Wakimoto module of $U_q(\widehat{sl}(M|N))$.

\section{Screening operator}
\label{Section:4}

In this Section 
we give the screening operators $Q_i~(1\leq i \leq M+N-1)$
that commute with $U_q(\widehat{sl}(M|N))$ for the level $c=k \neq -g$.
We define bosonic operators $S_i(z)~(1\leq i \leq M+N-1)$ that we call the screening currents as follows.
\begin{eqnarray}
S_i(z)&=&\sum_{j=i+1}^M \frac{e_{i,j}}{(q-q^{-1})z}(S_{i,j}^-(z)-S_{i,j}^+(z))+\sum_{j=M+1}^{M+N}e_{i,j}S_{i,j}(z)~~~~~(1\leq i \leq M-1),
\label{def:S1}\\
S_M(z)&=&\sum_{j=M+1}^{M+N}e_{M,j} S_{M,j}(z),
\label{def:S2}\\
S_i(z)&=&\sum_{j=i+1}^{M+N}\frac{e_{i,j}}{(q-q^{-1})z}(S_{i,j}^-(z)-S_{i,j}^+(z))~~~~~(M+1\leq i \leq M+N-1),
\label{def:S3}
\end{eqnarray}
where we set 
\begin{eqnarray}
S_{i,j}^\pm(z)=:e^{-(\frac{1}{k+g}a^i)(z;\frac{k+g}{2})}\widetilde{S}_{i,j}^\pm(z):,~~~
S_{i,j}(z)=:e^{-(\frac{1}{k+g}a^i)(z;\frac{k+g}{2})}\widetilde{S}_{i,j}(z):.
\end{eqnarray}
Here we set $e_{i,j}$ as follows.
\begin{eqnarray}
e_{i,i+1}&=&\left\{\begin{array}{cc}
1/d_{i,i}^2& (1\leq i \leq M-1),\\
-q^{-N+1}/d_{M,M}^2& (i=M),\\
-1/d_{i,i}^2& (M+1\leq i \leq M+N-1),
\end{array}\right.
\label{def:e1}
\\
e_{i,j}&=&
\left\{\begin{array}{cc}
1/d_{i,j}^3& (1\leq i \leq M-1, i+2 \leq j \leq M),\\
q^{k+1+M-N}/d_{i,j}^3& (1\leq i \leq M-1, M+1\leq j \leq M+N),\\
-q^{j-M-N}/d_{M,j}^3& (i=M, M+2\leq j \leq M+N),\\
-1/d_{i,j}^3& (M+1\leq i \leq M+N-1, i+2\leq j \leq M+N).
\label{def:e2}
\end{array}
\right.
\end{eqnarray}
For $1\leq i \leq M-1$ and $i+1\leq j\leq M$ we set
\begin{eqnarray}
\tilde{S}_{i,j}^\pm(z)&=&:e^{(b+c)^{i+1,j}(q^{M-N-j}z)-b_\pm^{i,j}(q^{M-N-j}z)
-(b+c)^{i,j}(q^{M-N-j \mp 1}z)}
\nonumber\\
&\times&
e^{\sum_{l=j+1}^{M}(\Delta_L^+ b_-^{i,l})(q^{M-N-l}z)+\sum_{l=M+1}^{M+N}(\Delta_L^+ b_-^{i,l})(q^{-M-N+l-1}z)}:.
\end{eqnarray}
For $1\leq i \leq M-1$ and $M+1\leq j\leq M+N$ we set
\begin{eqnarray}
\tilde{S}_{i,j}(z)&=&:e^{b^{i,j}(q^{-M-N+j}z)+b_+^{i+1,j}(q^{-M-N+j}z)-b^{i+1,j}(q^{-M-N+j+1}z)}\nonumber\\
&\times&e^{\sum_{l=j+1}^{M+N}(\Delta_L^+ b_-^{i,l})(q^{-M-N-1+l}z)}:.
\end{eqnarray}
For $M+1\leq j \leq M+N$ we set
\begin{eqnarray}
\tilde{S}_{M,j}(z)=:e^{(b+c)^{M+1,j}(q^{-M-N+j}z)+b^{M,j}(q^{-M-N+j}z)-\sum_{l=j+1}^{M+N}
(\Delta_L^0 b_-^{M,l})(q^{-M-N-1+l}z)}:.
\end{eqnarray}
For $M+1 \leq i \leq M+N-1$ and $i+1 \leq j \leq M+N$ we set
\begin{eqnarray}
\tilde{S}_{i,j}^\pm(z)&=&:e^{(b+c)^{i+1,j}(q^{-M-N+j}z)+b_\pm^{i,j}(q^{-M-N+j}z)-(b+c)^{i,j}(q^{-M-N+j\pm 1}z)}\nonumber\\
&\times&
e^{-\sum_{l=j+1}^{M+N}(\Delta_L^- b_-^{i,l})(q^{-M-N+l}z)}:.
\end{eqnarray}
The ${\bf Z}_2$-grading of the screening current is given by
$p(S_{M,j}(z))\equiv 1 \pmod{2}$ for $M+1 \leq j \leq M+N$ and zero otherwise.
The following is {\bf the
second main result} of this article.

\begin{thm}~~~\label{thm:2}
The screening currents $S_i(z)~(1\leq i \leq M+N-1)$ defined in (\ref{def:S1}), (\ref{def:S2}), and (\ref{def:S3})
commute with $U_q(\widehat{sl}(M|N))$ up to total difference.
\begin{eqnarray}
~[S_i(z),H^j]&=&0,\\
~[S_i(z),H_m^j]&=&0~~~(m \in {\bf Z}),\\
~[S_i(z_1),X^{+,j}(z_2)]&=&0,\\
~[S_i(z_1),X^{-,j}(z_2)]&=&\frac{\delta_{i,j}}{(q-q^{-1})z_1z_2}(\delta(q^{k+g}z_2/z_1)-\delta(q^{-k-g}z_2/z_1))\nonumber\\
&\times&:e^{-(\frac{a^i}{k+g})(z_1|-\frac{k+g}{2})}:.
\end{eqnarray}
\end{thm}
These screening currents reproduces those of $U_q(\widehat{sl}(M|1))$ 
upon the specialization $N=1$ \cite{ZG, Koj3}.

The $q$-difference operator with a parameter $\alpha$ is defined by
\begin{eqnarray}
{_\alpha \partial_z f(z)}=\frac{f(q^\alpha z)-f(q^{-\alpha}z)}{(q-q^{-1})z}.
\end{eqnarray}
The Jackson integral with parameters $p \in {\bf C}~(|p|<1)$ and $s \in {\bf C}^*$ is defined by
\begin{eqnarray}
\int_0^{s \infty}f(w)d_p w=s(1-p)\sum_{n \in {\bf Z}}f(sp^n) p^n.
\end{eqnarray}
The Jackson integral of the $q$-difference satisfy the following property.
\begin{eqnarray}
{\int_0^{s \infty}}{_\alpha \partial_w} f(w) d_pw=0~~~(p=q^{2\alpha}).
\label{integral-0}
\end{eqnarray}
We define the screening operators $Q_i~(1\leq i \leq M+N-1)$ as follows,
when the Jackson integrals are convergent.
\begin{eqnarray}
Q_i=\int_0^{s \infty}S_i(w)d_{q^{2(k+g)}}w.\label{sef:screening}
\end{eqnarray} 

\begin{cor}~~~The screening operators $Q_i$ $(1\leq i \leq M+N-1)$ commute with $U_q(\widehat{sl}(M|N))$.
\begin{eqnarray}
[Q_i,U_q(\widehat{sl}(M|N))]=0.
\end{eqnarray}
\end{cor}

For $r>0$ we define the theta function $[u]_r$ as
\begin{eqnarray}
[u]_r=q^{\frac{u^2}{r}-u}\frac{ 
\Theta_{q^{2r}}(q^{2u})}{(q^{2r};q^{2r})_\infty},
\end{eqnarray}
where we set
\begin{eqnarray}
\Theta_p(z)=(p;p)_\infty (z;p)_\infty (p/z;p)_\infty,~~~(z;p)_\infty=\prod_{m=0}^\infty (1-p^m z).
\end{eqnarray}
The theta function $[u]_r$ satisfies the following
quasi-periodicity condition
\begin{eqnarray}
[u+r]_r=-[u]_r,~~
[u+r\tau]_r=-e^{-\pi \sqrt{-1}\tau-\frac{2\pi\sqrt{-1}}{r}u}[u]_r,
\end{eqnarray}
where $\tau \in {\bf C}$ is given by $q=e^{-\pi \sqrt{-1}{/r \tau}}$.

\begin{prop}~~~
The screening currents $S_i(z)~(1\leq i \leq M+N-1)$ satisfy
\begin{eqnarray}
\left[u_1-u_2+\frac{A_{i,j}}{2}\right]_{k+g}S_i(z_1)S_j(z_2)=
\left[u_2-u_1+\frac{A_{i,j}}{2}\right]_{k+g}S_j(z_2)S_i(z_1),
\end{eqnarray}
where $q^{2u_j}=z_j~(j=1,2)$.
\end{prop}

\section{Proof of main results}
\label{Section:5}

In this Section we will show Theorem \ref{thm:1} and Theorem \ref{thm:2}.

\subsection{Proof of Theorem \ref{thm:1}}

We will show
\begin{eqnarray}
[X^{+,i}(z_1),X^{-,i}(z_2)]=\frac{1}{(q-q^{-1})z_1z_2}
\left(\delta(q^kz_2/z_1)\Psi_+^i(q^{\frac{k}{2}}z_2)-
\delta(q^{-k}z_2/z_1)\Psi_-^i(q^{-\frac{k}{2}}z_2)\right)
\label{eqn:1}
\end{eqnarray}
for $1\leq i \leq M+N-1$.\\
$\bullet$~For $1\leq i \leq M-1$ we have
\begin{eqnarray}
&&[E_{i,i}^+(z_1),F_{i,i}^{2,+}(z_2)]=-(q-q^{-1})\delta(q^k z_2/z_1)
:E_{i,i}^+(z_1)F_{i,i}^{2,+}(z_2):,\\
&&[E_{i,j}^\pm(z_1),F_{i,j}^{1,\pm}(z_2)]=\mp (q-q^{-1})\delta(q^{-k-2j+1\mp 1}z_2/z_1)
:E_{i,j}^\pm(z_1)F_{i,j}^{1,\pm}(z_2):~(1\leq j \leq i-1),\\
&&[E_{i,i}^-(z_1),F_{i,i}^{2,-}(z_2)]=(q-q^{-1})\delta(q^{-k-2i+2}z_2/z_1):
E_{i,i}^-(z_1)F_{i,i}^{2,-}(z_2):.
\end{eqnarray}
The remaining commutators vanish.
Hence we have
\begin{eqnarray}
&&[X^{+,i}(z_1),X^{-,i}(z_2)]\times (q-q^{-1})z_1z_2\nonumber\\
&=&\sum_{j=1}^{i-2}
\delta(q^{-k-2j}z_2/z_1)\left(c_{i,j}d_{i,j}^1:E_{i,j}^+(z_1)F_{i,j}^{1,+}(z_2):-c_{i,j+1}d_{i,j+1}^1 
:E_{i,j+1}^-(z_1)F_{i,j+1}^{1,-}(z_2):\right)\nonumber
\\
&+&\delta(q^{-k-2i+2}z_2/z_1)\left(
c_{i,i-1}d_{i,i-1}^1:E_{i,i-1}^+(z_1)F_{i,j+1}^{1,-}(z_2):
-c_{i,i}d_{i,i}^2 :E_{i,i}^-(z_1)f_{i,i}^{2,-}(z_2):
\right)\nonumber\\
&+&\delta(q^kz_2/z_1)
c_{i,i}d_{i,i}^2:E_{i,i}^+(z_1) F_{i,i}^{2,+}(z_2):
-c_{i,1}d_{i,1}^1 \delta(q^{-k}z_2/z_1):E_{i,1}^-(z_1)F_{i,1}^{1,-}(z_2):.
\end{eqnarray}
For $1\leq i \leq M-1$ we have
\begin{eqnarray}
&&:E_{i,j}^+(q^{-k-2j}z)F_{i,j}^{1,+}(z):=:E_{i,j+1}^-(q^{-k-2j}z)F_{i,j+1}^{1,-}(z):~~(1\leq j \leq i-2),\\
&&:E_{i,i-1}^+(q^{-k-2i+2}z)F_{i,i-1}^{1,+}(z):=
:E_{i,i}^-(q^{-k-2i+2}z)F_{i,i}^{2,-}(z):,\\
&&:E_{i,i}^+(q^kz)F_{i,i}^{2,+}(z):=\Psi_+^i(q^{\frac{k}{2}}z),~~~
:E_{i,1}^-(q^{-k}z)F_{i,1}^{1,-}(z):=\Psi_-^i(q^{-\frac{k}{2}}z),
\end{eqnarray}
and
\begin{eqnarray}
&&c_{i,j}d_{i,j}^1=c_{i,j+1}d_{i,j+1}^1~~~(1\leq j \leq i-2),\\
&&c_{i,i-1}d_{i,i-1}^1=c_{i,i}d_{i,i}^2,\\
&&c_{i,i}d_{i,i}^2=c_{i,1}d_{i,1}^1=1.
\end{eqnarray}
Hence we have (\ref{eqn:1}) for $1\leq i \leq M-1$.
\\
$\bullet$~For $i=M$ we have
\begin{eqnarray}
&&[E_{M,M}(z_1),F_{M,M}^{2,+}(z_2)]=\frac{1}{q^{M-1}z_1}\delta(q^k z_2/z_1):
E_{M,M}(z_1)F_{M,M}^{2,+}(z_2):,\\
&&[E_{M,j}(z_1),F_{M,j}^{1,\pm}(z_2)]=\frac{1}{q^{j-1}z_1}\delta(q^{-k-2j+1\mp 1}z_2/z_1):
E_{M,j}(z_1)F_{M,j}^{1,\pm}(z_2):~(1\leq j \leq M-1),
\\
&&[E_{M,M}(z_1),F_{M,M}^{2,-}(z_2)]=\frac{1}{q^{M-1}z_1}\delta(q^{-k-2M+2}z_2/z_1):
E_{M,M}(z_1)F_{M,M}^{2,-}(z_2):.
\end{eqnarray}
The remaining anti-commutators vanish.
Hence we have
\begin{eqnarray}
&&[X^{+,M}(z_1),X^{-,M}(z_2)]\times (q-q^{-1})z_1z_2\nonumber\\
&=&\sum_{j=1}^{M-2}\delta(q^{-k-2j}z_2/z_1)\left(
\frac{c_{M,j+1}d_{M,j+1}^1}{q^j}:E_{M,j+1}(z_1)F_{M,j+1}^{1,-}(z_2):
-\frac{c_{M,j}d_{M,j}^1}{q^{j-1}}:E_{M,j}(z_1)F_{M,j}^{1,+}(z_2):\right)\nonumber\\
&+&
\delta(q^{-k-2M+2}z_2/z_1)\nonumber\\
&\times&
\left(
\frac{c_{M,M}d_{M,M}^2}{q^{M-1}}:E_{M,M}(z_1)F_{M,M}^{2,-}(z_2):
-\frac{c_{M,M-1}d_{M,M-1}^1}{q^{M-2}}:E_{M,M-1}(z_1)F_{M,M-1}^{1,+}(z_2):\right)
\nonumber\\
&-&\delta(q^kz_2/z_1)\frac{c_{M,M}d_{M,M}^2}{q^{M-1}}:E_{M,M}(z_1)F_{M,M}^{2,+}(z_2):
+\delta(q^{-k}z_2/z_1)c_{M,1}d_{M,1}^1 :E_{M,1}(z_1)F_{M,1}^{1,-}(z_2):.
\end{eqnarray}
Moreover we have
\begin{eqnarray}
&&
:E_{M,j}(q^{-k-2j}z)F_{M,j}^{1,+}(z):=:E_{M,j+1}(q^{-k-2j}z)F_{M,j+1}^{1,-}(z):~~(1\leq j \leq M-2),\\
&&
:E_{M,M-1}(q^{-k-2M+2}z)F_{M,M-1}^{1,+}(z):
=:E_{M,M}(q^{-k-2M+2}z)F_{M,M}^{2,-}(z):,\\
&&
:E_{M,M}(q^{k}z)F_{M,M}^{2,+}(z):=\Psi_+^M(q^{\frac{k}{2}}z),~~~
:E_{M,1}(q^{-k}z)F_{M,1}^{1,-}(z):=\Psi_-^M(q^{-\frac{k}{2}}z),
\end{eqnarray}
and
\begin{eqnarray}
&&q c_{M,j}d_{M,j}^1=c_{M,j+1}d_{M,j+1}^1~~~(1\leq j \leq M-2),\\
&&q c_{M,M-1}d_{M,M-1}^1=c_{M,M}d_{M,M}^2,\\
&&c_{M,M}d_{M,M}^2=-q^{M-1},~~~
c_{M,1}d_{M,1}^1=-1.
\end{eqnarray}
Hence we have (\ref{eqn:1}) for $i=M$.
\\
$\bullet$~For $M+1 \leq i \leq M+N-1$ we have
\begin{eqnarray}
~[E_{i,i}^+(z_1), F_{i,i}^{2,+}(z_2)]
&=&(q-q^{-1})\delta(q^kz_2/z_1)
:E_{i,i}^+(z_1)F_{i,i}^{2,+}(z_2):,\\
~[E_{i,j}(z_1), F_{i,j}^1(z_2)]
&=&\frac{q^{k+1}}{(q-q^{-1})z_1z_2}(\delta(q^{-k-2j+2}z_2/z_1)-\delta(q^{-k-2j}z_2/z_1))\nonumber\\
&\times&
:E_{i,j}(z_1)F_{i,j}^1(z_2):~~~(1\leq j \leq M),
\\
~[E_{i,j}^\pm(z_1), F_{i,j}^{1,\pm}(z_2)]
&=&\pm (q-q^{-1})\delta(q^{-k-4M+2j-1\pm 1}z_2/z_1)
:E_{i,j}^\pm(z_1)F_{i,j}^{1,\pm}(z_2):\nonumber\\
&&~~~(M+1\leq j \leq i-1),\\
~[E_{i,i}^-(z_1), F_{i,i}^{2,-}(z_2)]
&=&-(q-q^{-1})\delta(q^{-k-4M+2i-2}z_2/z_1)
:E_{i,i}^-(z_1)F_{i,i}^{2,-}(z_2):.
\end{eqnarray}
The remaining commutators vanish.
Hence we have
\begin{eqnarray}
&&[X^{+,i}(z_1),X^{-,i}(z_2)]\times (q-q^{-1})z_1z_2\nonumber\\
&=&\sum_{j=1}^{M-1}
\delta(q^{-k-2j}z_2/z_1)q^{k+1}\left(c_{i,j+1}d_{i,j+1}^1 :E_{i,j+1}(z_1)F_{i,j+1}^1(z_2):-
c_{i,j}d_{i,j}^1:E_{i,j}(z_1)F_{i,j}^1(z_2):
\right)\nonumber\\
&+&\delta(q^{-k-2M}z_2/z_1)\left(
c_{i,M+1}d_{i,M+1}^1:E_{i,M+1}^-(z_1)F_{i,M+1}^{1,-}(z_2):
-c_{i,M}d_{i,M}^1q^{k+1}:E_{i,M}(z_1)F_{i,M}^1(z_2):
\right)\nonumber\\
&+&
\sum_{j=M+1}^{i-2}
\delta(q^{-k-4M+2j}z_2/z_1)\left(
c_{i,j+1}d_{i,j+1}^1 :E_{i,j+1}^-(z_1)F_{i,j+1}^{1,-}(z_2):
-c_{i,j}d_{i,j}^1 :E_{i,j}^+(z_1)F_{i,j}^{1,+}(z_2):
\right)\nonumber\\
&+&
\delta(q^{-k-4M+2i-2}z_2/z_1)\left(
c_{i,i}d_{i,i}^2 :E_{i,i}^-(z_1)F_{i,i}^{2,-}(z_2):
-c_{i,i-1}d_{i,i-1}^1 :E_{i,i-1}^+(z_1)F_{i,i-1}^1(z_2):\right)\nonumber\\
&-&\delta(q^kz_2/z_1) c_{i,i}d_{i,i}^2 :E_{i,i}^+(z_1)F_{i,i}^{2,+}(z_2):
+c_{i,1}d_{i,1}^1 q^{k+1} \delta(q^{-k}z_2/z_1):E_{i,1}(z_1)F_{i,1}^1(z_2):.
\end{eqnarray}
For $M+1\leq i \leq M+N-1$ we have
\begin{eqnarray}
&&
:E_{i,j+1}(q^{-k-2j}z)F_{i,j+1}(z):=
:E_{i,j}(q^{-k-2j}z)F_{i,j}^1(z):~~~(1\leq j \leq M-1),\\
&&
:E_{i,M+1}^-(q^{-k-2M}z)F_{i,M+1}^{1,-}(z):=:E_{i,M}(q^{-k-2M}z)F_{i,M}^1(z):,\\
&&
:E_{i,j+1}^-(q^{-k-4M+2j}z)F_{i,j+1}^{1,-}(z):=:E_{i,j}^+(q^{-k-4M+2j}z)F_{i,j}^{1,+}(z):~~~(M+1 \leq j \leq i-2),\\
&&
:E_{i,i}^-(q^{-k-4M+2i-2}z)F_{i,i}^{2,-}(z):=:E_{i,i-1}^+(q^{-k-4M+2i-2}z)F_{i,i-1}^{1,+}(z):,\\
&&
:E_{i,i}^+(q^kz)F_{i,i}^{2,+}(z):=\Psi_+^i(q^{\frac{k}{2}}z),~~~
:E_{i,1}(q^{-k}z)F_{i,1}^1(z):=\Psi_-^i(q^{-\frac{k}{2}}z),
\end{eqnarray}
and
\begin{eqnarray}
&&c_{i,j+1}d_{i,j+1}^1=c_{i,j}d_{i,j}^1~~~(1\leq j \leq M-1~{\rm or}~M+1\leq i-2),\\
&&c_{i,M+1}d_{i,M+1}^1=q^{k+1} c_{i,M}d_{i,M},\\
&&c_{i,i}d_{i,i}^2=c_{i,i-1}d_{i,i-1}^1,\\
&&c_{i,i}d_{i,i}^2=q^{k+1}c_{i,1}d_{i,1}^1=-1. 
\end{eqnarray}
Hence we have (\ref{eqn:1}) for $M+1\leq i \leq M+N-1$.
Now we have shown (\ref{eqn:1}) for all $1\leq i \leq M+N-1$.
Other defining relations of the Drinfeld realization
(\ref{def:Drinfeld1})-(\ref{def:Drinfeld10})
are shown in the same way.
We summarize useful formulae for proof of theorem in Appendix \ref{appendix:A}.

\subsection{Proof of Theorem \ref{thm:2}}

First, we will show
\begin{eqnarray}
[S_i(z_1), X^{-,i}(z_2)]=\frac{1}{(q-q^{-1})z_1z_2}\left(
\delta(q^{k+g}z_2/z_1)-\delta(q^{-k-g}z_2/z_1)
\right):e^{-\left(\frac{1}{k+g}a^i\right)\left(z_1|-\frac{k+g}{2}\right)}:
\label{eqn:2}
\end{eqnarray}
for $1\leq i \leq M+N-1$.\\
$\bullet$~For $1 \leq i \leq M-1$ we have
\begin{eqnarray}
~[S_{i,i+1}^+(z_1),F_{i,i}^{2,+}(z_2)]&=&(q-q^{-1})
\delta(q^{k-M+N+2i+2}z_2/z_1):
S_{i,i+1}^+(z_1)F_{i,i}^{2,+}(z_2):,
\\
~[S_{i,i+1}^-(z_1),F_{i,i}^{2,-}(z_2)]&=&-(q-q^{-1})
\delta(q^{k-M+N}z_2/z_1):
S_{i,i+1}^-(z_1)F_{i,i}^{2,-}(z_2):,
\\
~[S_{i,j}(z_1),F_{i,j}^3(z_2)]
&=&
\frac{q^{-k-M+N-1}}{(q-q^{-1})z_1z_2}
\left(\delta(q^{k+3M+N-2j}z_2/z_1)-\delta(q^{k+3M+N+2-2j}z_2/z_1)\right)\nonumber\\
&\times& :S_{i,j}(z_1)F_{i,j}^3(z_2):~~~~~(M+1 \leq j \leq M+N).
\end{eqnarray}
For $1\leq i \leq M-2$ and $i+2 \leq j \leq M$ we have
\begin{eqnarray}
&&[S_{i,j}^\pm(z_1),F_{i,j}^{3,\pm}(z_2)]=\pm (q-q^{-1})\delta(q^{k-M+N+2j-1\pm 1}z_2/z_1)
:S_{i,j}^\pm(z_1)F_{i,j}^{3,\pm}(z_2):.
\end{eqnarray}
The remaining commutators vanish.
Hence we have
\begin{eqnarray}
&&[S_i(z_1),X^{-,i}(z_2)]\times (q-q^{-1})z_1z_2\nonumber\\
&=&\sum_{j=i+2}^{M-1}\delta(q^{k-M+N+2j}z_2/z_1)\left(
e_{i,j}d_{i,j}^3:S_{i,j}^+(z_1)F_{i,j}^{3,+}(z_2):-e_{i,j+1}d_{i,j+1}^3
:S_{i,j+1}^-(z_1)F_{i,j+1}^{3,-}(z_2):
\right)\nonumber\\
&+&\delta(q^{k-M+N+2i+2}z_2/z_1)\left(
e_{i,i+1}d_{i,i}^2 :S_{i,i+1}^+(z_1)F_{i,i}^{2,+}(z_2):
-e_{i,i+2}d_{i,i+2}^3
:S_{i,i+2}^-(z_1)F_{i,i+2}^{3,-}(z_2):
\right)\nonumber\\
&+&\delta(q^{k+M+N}z_2/z_1)\nonumber\\
&\times&\left(
e_{i,M}d_{i,M}^3 :S_{i,M}^+(z_1)F_{i,M}^{3,+}(z_2):
-q^{-k-M+N-1}e_{i,M+1}d_{i,M+1}^3
:S_{i,M+1}(z_1)F_{i,M+1}^3(z_2):
\right)\nonumber\\
&+&\sum_{j=M+1}^{M+N-1}
\delta(q^{k+3M+N-2j}z_2/z_1)q^{-k-M+N-1}\nonumber\\
&\times&
\left(e_{i,j+1}d_{i,j+1}^3 :S_{i,j+1}(z_1)F_{i,j+1}^3(z_2):
-e_{i,j}d_{i,j}^3:S_{i,j}(z_1)F_{i,j}^3(z_2):\right)\nonumber\\
&+&\delta(q^{k+g}z_2/z_1)e_{i,M+N}d_{i,M+N}^3 q^{-k-M+N-1}:S_{i,M+N}(z_1)F_{i,M+N}^{3}(z):
\nonumber\\
&-&
\delta(q^{-k-g}z_2/z_1)e_{i,i+1}d_{i,i}^2 :S_{i,i+1}(z_1)F_{i,i}^{2,-}(z_2):.
\end{eqnarray}
Moreover we have
\begin{eqnarray}
:S_{i,i+1}^-(q^{-k-M+N}z)F_{i,i}^{2,-}(z):
&=&:S_{i,M+N}(q^{k+M-N}z)F_{i,M+N}^{3}(z):\nonumber\\
&=&:e^{-(\frac{1}{k+g}a^i)(z|-\frac{k+g}{2})}:~~~(1\leq i \leq M-1),
\\
:S_{i,j}^+(q^{k-M+N+2j}z)F_{i,j}^{3,+}(z):
&=&:S_{i,j+1}^-(q^{k-M+N+2j}z)F_{i,j+1}^{3,+}(z):\nonumber\\
&&~~~(1\leq i \leq M-2, i+2\leq j \leq M-1),
\\
:S_{i,i+1}^+(q^{k-M+N+2i+2}z)F_{i,i}^{2,+}(z):
&=&:S_{i,i+2}^-(q^{k-M+N+2i+2}z)F_{i,i+2}^{3,-}(z):~~(1 \leq i \leq M-2),\\
:S_{i,M}^+(q^{k+M+N}z)F_{i,M}^{3,+}(z):
&=&:S_{i,M+1}(q^{k+M+N}z)F_{i,M+1}^3(z):~~(1\leq i \leq M-2),\\
:S_{i,j+1}(q^{k+3M+N-2j}z)F_{i,j+1}^3(z):
&=&:S_{i,j}(q^{k+3M+N-2j}z)F_{i,j}^{3}(z):\nonumber\\
&&~~(1\leq i \leq M-1, M+1\leq j \leq M+N-1),
\end{eqnarray}
and
\begin{eqnarray}
&&
e_{i,i+1}d_{i,i}^2=e_{i,i+2}d_{i,i+2}^3~~(1\leq i \leq M-1),\\
&&
e_{i,j}d_{i,j}^3=e_{i,j+1}d_{i,j+1}^3~~\left(\begin{array}{c}
1\leq i \leq M-2, i+2 \leq j \leq M+N-1\\
{\rm or}~~
1\leq i \leq M-1, M+1\leq j \leq M+N-1 \end{array}
\right),
\\
&&
e_{i,M}d_{i,M}^3=e_{i,M+1}d_{i,M+1}^3 q^{-k-M+N-1},~~~(1 \leq i \leq M-1),
\\
&&
e_{i,M+N}d_{i,M+N}^3q^{-k-M+N-1}=e_{i,i+1}d_{i,i}^2=1~~~(1\leq i \leq M-1).
\end{eqnarray}
Hence we have (\ref{eqn:2}) for $1\leq i \leq M-1$.
\\
$\bullet$~For $i=M$ and $M+2\leq j \leq M+N$ we have
\begin{eqnarray}
&&[S_{M,M+1}(z_1),F_{M,M}^{2,\pm}(z_2)]=\frac{q^{N-1}}{z_1}\delta(q^{k+M+N-1 \mp 1}z_2/z_1):
S_{M,M+1}(z_1)F_{M,M}^{2,\pm}(z_2):,
\\
&&[S_{M,j}(z_1),F_{M,j}^{3,\pm}(z_2)]=\frac{q^{M+N-j}}{z_1}\delta(q^{k+3M+N-2j+1 \mp 1}z_2/z_1):
S_{M,j}(z_1)F_{M,j}^{3,\pm}(z_2):.
\end{eqnarray}
The remaining anti-commutators vanish.
Hence we have
\begin{eqnarray}
&&[S_M(z_1),X^{-,M}(z_2)]\times (q-q^{-1})z_1z_2\nonumber\\
&=&\delta(q^{k+M+N-2}z_2/z_1)\nonumber\\
&\times&\left(
q^{N-2} e_{M,M+2}d_{M,M+2}^3
:S_{M,M+1}(z_1)F_{M,M}^{2,+}(z_2):-q^{N-1}e_{M,M+1}d_{M,M}^2
:S_{M,M+1}(z_1)F_{M,M}^{2,+}(z_2):
\right)\nonumber\\
&+&
\sum_{j=M+2}^{M+N-1}\delta(q^{k+3M+N-2j}z_2/z_1)\nonumber\\
&\times&
\left(q^{M+N-j-1}e_{M,j+1}d_{M,j+1}^3 :S_{M,j+1}(z_1)F_{M,j+1}^{3,-}(z_2):-q^{M+N-j}
e_{M,j}d_{M,j}^3 :S_{M,j}(z_1)F_{M,j}^{3,+}(z_2):
\right)\nonumber\\
&-&\delta(q^{k+g}z_2/z_1)e_{M,M+N}d_{M,M+N}^3 :S_{M,M+N}(z_1)F_{M,M+N}^{3,+}(z_2):\nonumber\\
&+&e_{M,M+1}d_{M,M}^2 \delta(q^{-k-g}z_2/z_1):S_{M,M+1}(z_1)F_{M,M}^2(z_2):.
\end{eqnarray}
For $M+2 \leq i \leq M+N-1$ we have
\begin{eqnarray}
:S_{M,M+1}(q^{-k-M+N}z)F_{M,M}^{2,-}(z):
&=&
:S_{M,M+N}(q^{k+M-N}z)F_{M,M+N}^{3,+}(z):\nonumber\\
&=&:e^{-(\frac{1}{k+g}a^M)(z|-\frac{k+g}{2})}:,
\\
:S_{M,j}(q^{k+3M+N-2j}z)F_{M,j}^{3,+}(z):
&=&:S_{M,j+1}(q^{k+3M+N-2j}z)F_{M,j+1}^{3,-}(z):,
\\
:S_{M,M+1}(q^{k+M+N-2}z)F_{M,M}^{2,+}(z):&=&
:S_{M,M+2}(q^{k+M-N-2}z)F_{M,M+2}^{3,-}(z):,
\end{eqnarray}
and
\begin{eqnarray}
&&e_{M,M+1} d_{M,M}^2q^{N-1}=-1,~~
e_{M,M+N}d_{M,M+N}^3=-1,\\
&&
qe_{M,M+1}d_{M,M}^2=e_{M,M+2}d_{M,M+2}^3,
\\
&&
q e_{M,j}d_{M,j}^3=e_{M,j+1}d_{M,j+1}^3.
\end{eqnarray}
Hence we have (\ref{eqn:2}) for $i=M$.\\
$\bullet$~For $M+1 \leq i \leq M+N-1$ and $i+2 \leq j \leq M+N$ we have
\begin{eqnarray}
&&
[S_{i,i+1}^+(z_1),F_{i,i}^{2,+}(z_2)]=-(q-q^{-1})\delta(q^{k+3M+N-2i-2}z_2/z_1)
:S_{i,i+1}^+(z_1)F_{i,i}^{2,+}(z_2):,\\
&&
[S_{i,i+1}^-(z_1),F_{i,i}^{2,-}(z_2)]=(q-q^{-1})\delta(q^{-k-M+N}z_2/z_1)
:S_{i,i+1}^-(z_1)F_{i,i}^{2,-}(z_2):,\\
&&
[S_{i,j}^\pm(z_1),F_{i,j}^{3,\pm}(z_2)]=\mp(q-q^{-1})\delta(q^{k+3M+N-2j+1\mp 1}z_2/z_1)
:S_{i,j}^\pm(z_1)F_{i,j}^{3,\pm}(z_2):.
\end{eqnarray}
The remaining commutators vanish.
Hence we have
\begin{eqnarray}
&&[S_i(z_1),X^{-,i}(z_2)]\times (q-q^{-1})z_1z_2\nonumber\\
&=&\delta(q^{k+3M+N-2i-2}z_2/z_1)\left(
e_{i,i+2}d_{i,i+2}^3:S_{i,i+2}^-(z_1)F_{i,i+2}^{3,-}(z_2):-
e_{i,i+1}d_{i,i}^2:S_{i,i+1}^+(z_1)F_{i,i}^{2,+}(z_2):
\right)\nonumber\\
&+&\sum_{j=i+2}^{M+N-1}\delta(q^{k+3M+N-2j}z_2/z_1)
\left(e_{i,j+1}d_{i,j+1}^3:S_{i,j+1}^-(z_1)F_{i,j+1}^{3,-}(z_2):
-e_{i,j}d_{i,j}^3:S_{i,j}^+(z_1)F_{i,j}^{3,+}(z_2):\right)\nonumber\\
&-&\delta(q^{k+g}z_2/z_1)e_{i,M+N}d_{i,M+N}^3 :S_{i,M+N}^+(z_1)F_{i,M+N}^{3,+}(z_2):\nonumber\\
&+&\delta(q^{-k-g}z_2/z_1)e_{i,i+1}d_{i,i}^2 :S_{i,i+1}^-(z_1)F_{i,i}^{2,-}(z_2):.
\end{eqnarray}
Moreover we have
\begin{eqnarray}
:S_{i,i+1}^-(q^{-k-M+N}z)F_{i,i}^{2,-}(z):&=&:S_{i,M+N}^+(q^{k+M-N}z)F_{i,M+N}^{3,+}(z):\nonumber\\
&=&:e^{-(\frac{1}{k+g}a^i)(z|-\frac{k+g}{2})}:~~~(M+1 \leq i \leq M+N-1),\\
:S_{i,j}^+(q^{k+3M+N-2j}z)F_{i,j}^{3,+}(z):&=&:S_{i,j+1}^-(q^{k+3M+N-2j}z)F_{i,j+1}^{3,-}(z):\nonumber\\
&&(M+1\leq i \leq M+N-1, i+2\leq j \leq M+N-1),
\\
:S_{i,i+2}^-(q^{k+3M+N-2i-2}z)F_{i,i+2}^{3,-}(z):&=&:S_{i,i+1}^+(q^{k+3M+N-2i-2}z)F_{i,i}^{2,+}(z):~~~(M+1 \leq i \leq M+N-2).\nonumber
\\
\end{eqnarray}
For $M+1 \leq i \leq M+N-1$ we have
\begin{eqnarray}
&&e_{i,i+1}d_{i,i}^2=-1,~~e_{i,M+N}d_{i,M+N}^3=-1,\\
&&e_{i,i+1}d_{i,i}^2=e_{i,i+2}d_{i.i+2}^3,~~e_{i,j}d_{i,j}^3=e_{i,j+1}d_{i,j+1}^3~~(i+3 \leq j \leq M+N-1).
\end{eqnarray}
Hence we have (\ref{eqn:2}) for $M+1\leq i \leq M+N-1$.
Now we have shown (\ref{eqn:2}) for all $1\leq i \leq M+N-1$.

~\\

Next, we will show
\begin{eqnarray}
[S_i(z_1), X^{+,i}(z_2)]=0~~~~~(1\leq i \leq M+N-1).
\label{eqn:3}
\end{eqnarray}
$\bullet$~For $1\leq i \leq M-1$ we have
\begin{eqnarray}
[S_{i,i+1}^\pm(z_1),E_{i,i}^\pm(z_2)]=\mp (q-q^{-1})\delta(q^{2i-M+N}z_2/z_1):
S_{i,i+1}^\pm(z_1)E_{i,i}^\pm(z_2):.
\end{eqnarray}
Hence we have
\begin{eqnarray}
[S_i(z_1), X^{+,i}(z_2)]=\delta(q^{2i-M+N}z_2/z_1)e_{i,i+1}c_{i,i}
\left(:S_{i,i+1}^+(z_1)E_{i,i}^+(z_2):
-:S_{i,i+1}^-(z_1)E_{i,i}^-(z_2):
\right).
\end{eqnarray}
Moreover we have
\begin{eqnarray}
:S_{i,i+1}^+(q^{2i-M+N}z)E_{i,i}^+(z):=
:S_{i,i+1}^-(q^{2i-M+N}z)E_{i,i}^-(z):.
\end{eqnarray}
Hence we have (\ref{eqn:3}) for $1\leq i \leq M-1$.\\
$\bullet$~For $i=M$ 
all anti-commutators vanish.
Hence we have (\ref{eqn:3}) for $i=M$.\\
$\bullet$~For $M+1 \leq i \leq M+N-1$ we have
\begin{eqnarray}
[S_{i,i+1}^\pm(z_1),E_{i,i}^\pm(z_2)]=\pm (q-q^{-1})\delta(q^{3M+N-2i}z_2/z_1):
S_{i,i+1}^\pm(z_1)E_{i,i}^\pm(z_2):.
\end{eqnarray}
Hence we have
\begin{eqnarray}
[S_i(z_1), X^{+,i}(z_2)]=
\delta(q^{3M+N-2i}z_2/z_1)e_{i,i+1}c_{i,i}
\left(:S_{i,i+1}^+(z_1)E_{i,i}^+(z_2):
-:S_{i,i+1}^-(z_1)E_{i,i}^-(z_2):
\right).
\end{eqnarray}
Moreover we have
\begin{eqnarray}
:S_{i,i+1}^+(q^{3M+N-2i}z)E_{i,i}^+(z):
=:S_{i,i+1}^-(q^{3M+N-2i}z)E_{i,i}^-(z):.
\end{eqnarray}
Hence we have (\ref{eqn:3}) for $M+1\leq i \leq M+N-1$.
Now we have shown (\ref{eqn:3}) for all $1\leq i \leq M+N-1$.

Other commutation relations of the screening currents
are shown in the same way.
We summarize useful formulae for proof of theorem in Appendix \ref{appendix:A}.

\section{Concluding remarks}
\label{Section:6}

In this article we found
a bosonization of $U_q(\widehat{sl}(M|N))$
for an arbitrary level $k \in {\bf C}$.
Our bosonization
is obtained
from a $q$-difference realization of $U_q(sl(M|N))$ 
(see Appendix \ref{appendix:B}) 
by the replacement
\begin{eqnarray}
q^{\vartheta_{i,j}} &\to& e^{\pm b_{\mp}^{i,j}(z)},
\\
x_{i,j} &\to& \left\{
\begin{array}{cc}
:e^{(b+c)^{i,j}(z)}: & (\nu_i \nu_j=+),\\
:e^{-b^{i,j}(z)}:~{\rm or}~:e^{-b_\pm^{i,j}(z)-b^{i,j}(z)}:& (\nu_i \nu_j=-),
\end{array}\right.
\\
q^{\lambda_i} &\to& e^{\pm a_{\pm}^i(z)},\\
~[\vartheta_{i,j}]_q &\to& 
\left\{\begin{array}{cc}
\frac{e^{\pm b_+^{i,j}(z)}-e^{\pm b_-^{i,j}(z)}}{(q-q^{-1})z} &(\nu_i \nu_j=+),\\
1 & (\nu_i\nu_j=-).
\end{array}\right. 
\end{eqnarray}
For instance, we have
\begin{eqnarray}
f_{M,j}^1
&\to&
\frac{1}{(q-q^{-1})z}:\left(e^{-b_-^{j,M}(z)-(b+c)^{j,M}(z)}-
e^{-b_-^{j,M}(z)-(b+c)^{j,M}(z)}\right)\nonumber
\\
&&\times e^{
a_-^M(z)-b_-^{j,M+1}(z)-b^{j,M+1}(z)
-\sum_{l=j+1}^{M-1}(\Delta_R^0 b_-^{l,M})(z)+\sum_{l=M+2}^{M+N}(\Delta_R^0 b_-^{M,l})(z)}:\nonumber
\\
&\to&
\frac{1}{(q-q^{-1})z}
:\left(
e^{-b_+^{j,M}(q^{-k-j}z)-(b+c)^{j,M}(q^{-k-j-1}z)}
-e^{-b_-^{j,M}(q^{-k-j}z)-(b+c)^{j,M}(q^{-k-j+1}z)}
\right)\nonumber\\
&&\times
e^{a_-^M(q^{-\frac{k+g}{2}}z)-b_-^{j,M+1}(q^{-k-j}z)-b^{j,M+1}(q^{-k-j+1}z)}
\nonumber\\
&&\times
e^{-\sum_{l=j+1}^{M-1}(\Delta_R^0 b_-^{l,M})(q^{-k-l}z)+\sum_{l=M+2}^{M+N}(\Delta_R^0 b_-^{M,l})(q^{-k-2M-1+l}z)}:\nonumber\\
&=&\frac{1}{(q-q^{-1})z}(F_{M,j}^{1,-}(z)-F_{M,j}^{1,+}(z))~~~~~(1\leq j \leq M-1).
\end{eqnarray}

Taking the limit $q \to 1$
we obtain a bosonization of the affine superalgebra $\widehat{sl}(M|N)$ 
for an arbitrary level $k \in {\bf C}$.
Bosonizations of the affine superalgebra $\widehat{sl}(M|N)$ for level $k$ have been studied in 
\cite{Bars, IR, BKT, BT, BFST, IK, DGZ, YZL1, YZL2}.
We compare our bosonization with those of \cite{YZL2}.
In the limit $q\to1$ we introduce operators $\alpha_i(z)~(1\leq i \leq M+N-1)$,
$\beta_{i,j}(z), \widehat{\beta}_{i,j}(z), \gamma_{i,j}(z)~(1\leq i<j \leq M+N, \nu_i\nu_j=+)$,
and $\psi_{i,j}(z), \psi_{i,j}^\dagger(z)~(1\leq i<j \leq M+N, \nu_i\nu_j=-)$ as follows.
\begin{eqnarray}
&&
\alpha_i(z)=\partial_z \left(a^i(z)\right),\\
&&
\beta_{i,j}(z)=:\partial_z \left(e^{-c^{i,j}(z)}\right) e^{-b^{i,j}(z)}:,~
\widehat{\beta}_{i,j}(z)=:\partial_z \left(e^{-b^{i,j}(z)}\right) e^{-c^{i,j}(z)}:,\\
&&
\gamma_{i,j}(z)=:e^{(b+c)^{i,j}(z)}:,
\\
&&
\psi_{i,j}(z)=:e^{b^{i,j}(z)}:,~~~\psi_{i,j}^\dagger (z)=:e^{-b^{i,j}(z)}:.
\end{eqnarray}
They satisfy the following relations.
\begin{eqnarray}
&&
\alpha_i(z)\alpha_j(w)=\frac{(k+g)A_{i,j}}{(z-w)^2}+\cdots,\\
&&
\beta_{i,j}(z)\gamma_{i',j'}(w)=\frac{\delta_{i,i'}\delta_{j,j'}}{z-w}+\cdots,
~~~
\gamma_{i,j}(z)\beta_{i',j'}(w)=-\frac{\delta_{i,i'}\delta_{j,j'}}{z-w}+\cdots,\\
&&
\widehat{\beta}_{i,j}(z)\gamma_{i',j'}(w)=-\frac{\delta_{i,i'}\delta_{j,j'}}{z-w}+\cdots,
~~~
\gamma_{i,j}(z)\widehat{\beta}_{i',j'}(w)=\frac{\delta_{i,i'}\delta_{j,j'}}{z-w}+\cdots,\\
&&
\psi_{i,j}(z)\psi_{i',j'}^\dagger (w)=\frac{\delta_{i,i'}\delta_{j,j'}}{z-w}+\cdots,
~~~
\psi_{i,j}^\dagger(z)\psi_{i',j'}(w)=\frac{\delta_{i,i'}\delta_{j,j'}}{z-w}+\cdots.
\end{eqnarray}
In the limit $q \to 1$ our bosonization becomes simpler
because the operators $a_\pm^i(z)$, $b_\pm^{i,j}(z)$, $(\Delta_L^\epsilon b_\pm^{i,j})(z)$ and
$(\Delta_R^\epsilon b_\pm^{i,j})(z)$ disappear.
In order to resolve a singularity in denominator $(q-q^{-1})$,
we take the limit of $\frac{1}{(q-q^{-1})z}(E_{i,j}^{+}(z)-E_{i,j}^{-}(z))$ instead of $\frac{1}{(q-q^{-1})z}E_{i,j}^{\pm}(z)$.
In Appendix \ref{appendix:C}
we summarize useful formulae to take the limit $q \to 1$ upon the condition $k \neq -g$.
Taking the limit $q\to 1$, we obtain the following bosonization of the affine superalgebra $\widehat{sl}(M|N)$.
In what follows we set the coefficients $c_{i,j}=1$ for simplicity.
\begin{eqnarray}
H^i(z)&=&\alpha_i(z)+\sum_{j=1}^i:(\widehat{\beta}_{j,i}(z)\gamma_{j,i}(z)-\widehat{\beta}_{j,i+1}(z)\gamma_{j,i+1}(z)):\nonumber\\
&&+\sum_{j=i+1}^M:(\widehat{\beta}_{i+1,j}(z)\gamma_{i+1,j}(z)-\widehat{\beta}_{i,j}(z)\gamma_{i,j}(z)):\nonumber\\
&&+\sum_{j=M+1}^{M+N}:( (\partial_z{\psi}_{i+1,j})(z) \psi_{i+1,j}^\dagger(z)-
(\partial_z {\psi}_{i,j})(z) \psi_{i,j}^\dagger(z)):~~~(1\leq i \leq M-1),
\\
H^M(z)&=&\alpha_M(z)+\sum_{j=1}^{M-1}
:((\partial_z {\psi}_{j,M+1})(z)\psi_{j,M+1}^\dagger(z)
+\widehat{\beta}_{j,M}(z)\gamma_{j,M}(z)):\nonumber\\
&&-\sum_{j=M+2}^{M+N}:(\widehat{\beta}_{M+1,j}(z)\gamma_{M+1,j}(z)+(\partial_z {\psi}_{M,j})(z) \psi_{M,j}^\dagger(z)):,
\\
H^i(z)&=&
\alpha_i(z)+\sum_{j=1}^M
:((\partial_z {\psi}_{j,i+1})(z)\psi_{j,i+1}^\dagger(z)-
(\partial_z{\psi}_{j,i})(z)\psi_{j,i}^\dagger(z)):\nonumber\\
&&+\sum_{j=M+1}^i:(\widehat{\beta}_{j,i+1}(z)\gamma_{j,i+1}(z)-\widehat{\beta}_{j,i}(z)\gamma_{j,i}(z)):\nonumber\\
&&+\sum_{j=i+1}^{M+N}
:(\widehat{\beta}_{i,j}(z)\gamma_{i,j}(z)-\widehat{\beta}_{i+1,j}(z)\gamma_{i+1,j}(z)):~~(M+1\leq i \leq M+N-1).
\end{eqnarray}
\begin{eqnarray}
X^{+,i}(z)&=&\sum_{j=1}^i :\beta_{j,i+1}(z)\gamma_{j,i}(z):~~~(1\leq i \leq M-1),\\
X^{+,M}(z)&=&\sum_{j=1}^M :\gamma_{j,M}(z)\psi_{j,M+1}(z):,\\
X^{+,i}(z)&=&\sum_{j=1}^M:\psi_{j,i+1}(z)\psi_{j,i}^\dagger(z):-\sum_{j=M+1}^i
:\beta_{j,i+1}(z)\gamma_{j,i}(z):~(M+1\leq i \leq M+N-1).
\end{eqnarray}
\begin{eqnarray}
X^{-,i}(z)&=&-:\alpha_i(z)\gamma_{i,i+1}(z):-\kappa_i :\partial_z \gamma_{i,i+1}(z):\nonumber\\
&&+\sum_{j=1}^{i-1}:\beta_{j,i}(z)\gamma_{j,i+1}(z):
-\sum_{j=i+2}^M:\beta_{i+1,j}(z)\gamma_{i,j}(z):
-\sum_{j=M+1}^{M+N}:\psi_{i+1,j}(z)\psi_{i,j}^\dagger(z):\nonumber\\
&&+\sum_{j=i+1}^M
:(\widehat{\beta}_{i,j}(z)\gamma_{i,j}(z)
-\widehat{\beta}_{i+1,j}(z)\gamma_{i+1,j}(z))\gamma_{i,i+1}(z):
\\
&&+\sum_{j=M+1}^{M+N}
:((\partial_z {\psi}_{i,j})(z)\psi_{i,j}^\dagger(z)-
(\partial_z{\psi}_{i+1,j})(z)\psi_{i+1,j}^\dagger(z))\gamma_{i,i+1}(z):~~~(1\leq i \leq M-1),
\nonumber\\
X^{-,M}(z)&=&:\alpha_M(z)\psi_{M,M+1}^\dagger(z):+\kappa_M :\partial_z \psi_{M,M+1}^\dagger(z):\nonumber\\
&&-\sum_{j=1}^{M-1}:\beta_{j,M}(z)\psi_{j,M+1}^\dagger(z):
-\sum_{j=M+2}^{M+N}:\beta_{M+1,j}(z)\psi_{M,j}^\dagger(z):\nonumber\\
&&-\sum_{j=M+2}^{M+N}
:(\widehat{\beta}_{M+1,j}(z)\gamma_{M+1,j}^\dagger(z)+
(\partial_z {\psi}_{M,j})(z)\psi_{M,j}^\dagger(z))\psi_{M,M+1}^\dagger(z):,\\
X^{-,i}(z)&=&
:\alpha_i(z)\gamma_{i,i+1}(z):+\kappa_i:\partial_z \gamma_{i,i+1}(z):\nonumber\\
&&-\sum_{j=1}^M:\psi_{j,i}(z)\psi_{j,i+1}^\dagger(z):
+\sum_{j=M+1}^{i-1}:\beta_{j,i}(z)\gamma_{j,i+1}(z):
-\sum_{j=i+2}^{M+N}:\beta_{i+1,j}(z)\gamma_{i,j}(z):\nonumber\\
&&+\sum_{j=i+1}^{M+N}:(\widehat{\beta}_{i,j}(z)\gamma_{i,j}(z)-
\widehat{\beta}_{i+1,j}(z)\gamma_{i+1,j}(z))\gamma_{i,j}(z):~~~(M+1\leq i \leq M+N-1).\nonumber\\
\end{eqnarray}
Here we have set the coefficients $\kappa_i$ by 
\begin{eqnarray}
\kappa_i=\left\{
\begin{array}{cc}
k+i& (1\leq i \leq M-1)\\
k+M-1& (i=M)\\
k+2M-i& (M+1\leq i \leq M+N-1)
\end{array}\right.
.
\end{eqnarray}
For $k \neq -g$ we have a bosonization of the screening current as follows.
\begin{eqnarray}
S_i(z)&=&\sum_{j=i+1}^M:\tilde{s}_i(z)\beta_{i,j}(z)\gamma_{i+1,j}(z):+
\sum_{j=M+1}^{M+N}:\tilde{s}_i(z)\psi_{i,j}(z)\psi_{i+1,j}^\dagger(z):~(1\leq i \leq M-1),\\
S_M(z)&=&\sum_{j=M+1}^{M+N}:\tilde{s}_M(z)\gamma_{M+1,j}(z)\psi_{M,j}(z):,\\
S_i(z)&=&\sum_{j=i+1}^{M+N}:\tilde{s}_i(z)\beta_{i,j}(z)\gamma_{i+1,j}(z):~~~(M+1\leq i \leq M+N-1).
\end{eqnarray}
Here we have set the operator $\tilde{s}_i(z)=:e^{-\left(\frac{1}{k+g}a^i\right)(z;0)}:$.
Our bosonization is similar as those of \cite{YZL2},
because both bosonizations are based on the same differential realization of $sl(M|N)$
(see Appendix \ref{appendix:B}).
The bosonization of
\cite{YZL2} is reproduced from our bosonization by the following "formal replacement"
\begin{eqnarray}
\widehat{\beta}_{i,j}(z) \to \beta_{i,j}(z),~~~
\partial_z {\psi}_{i,j}(z) \to \psi_{i,j}(z),~~~
\kappa_i \to 
\left\{
\begin{array}{cc}
k+i-1& (1\leq i \leq M-1)\\
k+M-1& (i=M)\\
k+M+1-i& (M+1\leq i \leq M+N-1)
\end{array}\right.,
\end{eqnarray}
with $\alpha_i(z), \beta_{i,j}(z), \gamma_{i,j}(z), \psi_{i,j}(z), \psi_{i,j}^\dagger(z)$ fixed.
Of course the map satisfying both $\partial_z \psi_{i,j}(z) \to \psi_{i,j}(z)$ and $\psi_{i,j}(z) \to \psi_{i,j}(z)$ is impossible.
This is the reason why we used the word "formal replacement".

In order to calculate correlation functions of exactly solvable models,
we have to prepare the vertex operator $\Phi_{V(\mu)}^{V(\nu) V_\lambda}(z)$
that satisfies intertwining property \cite{JM}.
We will propose a bosonization of the vertex operator.
In what follows we assume $k+g\neq 0$ and $g \neq 0$.
We define a bosonic operator
$\phi^{\bar{\lambda}}(z)$ 
for the weight $\bar{\lambda}=\sum_{i=1}^{M+N-1}l_i \bar{\Lambda}_i$
as follows.
\begin{eqnarray}
\phi^{\bar{\lambda}}(z)
=:\exp\left(\sum_{i,j=1}^{M+N-1}
\left(\frac{l_i}{k+g}\cdot \frac{\alpha_{i,j}}{g} \cdot \frac{\beta_{i,j}}{1}~a^j\right)\left(z;-\frac{k+g}{2}\right)\right):,
\end{eqnarray}
where we have set
\begin{eqnarray}
&&\alpha_{i,j}=\left\{\begin{array}{cc}
{\rm min}(i,j) & ({\rm min}(i,j)\leq M+1),\\
2(M+1)-{\rm min}(i,j) & ({\rm min}(i,j) \geq M+2),
\end{array}\right.\\
&&\beta_{i,j}=\left\{
\begin{array}{cc}
M-N-{\rm max}(i,j) & ({\rm max}(i,j) \leq M+1),\\
-M-N-2+{\rm max}(i,j) & ({\rm max}(i,j) \geq M+2).
\end{array}
\right.
\end{eqnarray}
The bosonic operator $\phi^{\bar{\lambda}}(z)$ satisfies the following relations for $1\leq i \leq M+N-1$.
\begin{eqnarray}
~[H_m^i, \phi^{\bar{\lambda}}(z)]
&=&\frac{1}{m}[l_i m]_q 
q^{-\frac{k}{2}|m|}z^m \phi^{\bar{\lambda}}(z)~~~(m \in {\bf Z}_{\neq 0}),\\
~[X^{+,i}(z), \phi^{\bar{\lambda}}(z)]&=&0,\\
(z_1-q^{l_i}z_2)X^{-,i}(z_1)\phi^{\bar{\lambda}}(z_2)&=&
(q^{l_i}z_1-z_2)\phi^{\bar{\lambda}}(z_2)X^{-,i}(z_1).
\end{eqnarray}
We set the bosonic operators $\phi^{\bar{\lambda}}_{i_1,i_2,\cdots,i_n}(z)~~(1\leq i_1, i_2, \cdots,i_n \leq M+N-1)$ as follows.
\begin{eqnarray}
\phi^{\bar{\lambda}}_{i_1,i_2,\cdots,i_n}(z)=\left[\phi^{\bar{\lambda}}_{i_1,i_2,\cdots,i_{n-1}}, X_0^{-,i_n}\right]_{q^x},
\label{def:bosonic1}
\end{eqnarray}
where $x=(\bar{\lambda}-\sum_{s=1}^{n-1} \bar{\alpha}_{i_s}|\bar{\alpha}_{i_n}))$.
The ${\bf Z}_2$-grading is given by 
$p(\phi^{\bar{\lambda}}(z))\equiv 0 \pmod{2}$ and 
$p(\phi^{\bar{\lambda}}_{i_1,i_2,\cdots,i_n}(z))\equiv \sum_{s=1}^n p(X_m^{-,i_s}) \pmod{2}$.
Let $\Phi_{V(\mu)}^{V(\nu) V_\lambda}(z)$ be the vertex operator
satisfying the intertwining property :
$\Phi_{V(\mu)}^{V(\nu) V_\lambda} : V(\mu) \to V(\nu) \otimes {V_\lambda}_z$.
Here $V(\mu)$ and $V(\nu)$ are highest weight modules.
$V_\lambda$ and
${V_{\lambda}}_z$
are a typical module and its evaluation module, respectively \cite{Palev-Tolstoy}.
Let $\Phi_{V(\mu)~i_1,i_2,\cdots,i_n}^{V(\nu)~V_\lambda}(z)$ be
$\Phi_{V(\mu)}^{V(\nu) V_\lambda}(z)=\sum_{i_1,i_2,\cdots,i_n}
\Phi_{V(\mu)~i_1,i_2,\cdots,i_n}^{V(\nu)~V_\lambda}(z)
v_{i_1,i_2,\cdots,i_n}$ where $\{v_{i_1,i_2,\cdots,i_n}\}$ is a basis of $V_{\lambda}$.
We propose a bosonization of the vertex operator as follows.
\begin{eqnarray}
\Phi_{V(\mu)~i_1,i_2,\cdots,i_n}^{V(\nu)~V_\lambda}(z)=
\eta_0 \xi_0 \prod_{j=1}^{M+N-1} : Q_j^{n_j} \cdot
\phi^{\bar{\lambda}}_{i_1,i_2,\cdots,i_n}(q^{k+g}z) : \eta_0 \xi_0,
\label{def:bosonic2}
\end{eqnarray}
where $n_j \in {\bf N}~(1\leq j \leq M+N-1)$.
Here $\eta_0 \xi_0$ is the projection operator on the Wakimoto module (see Section \ref{Section:3})
and $Q_j$ are the screening operators (see Section \ref{Section:4}).
In order to balance "background charge" of the Wakimoto module, we multiply the screening operators.
Trace of the vertex operators vanishes if we do not multiply the screening operators.
For small rank we have checked this conjecture by direct calculation.
For $2 \leq M  \leq 4$ and $N=1$,
the intertwining property of the vertex operator was checked \cite{ZG, Koj3}.
We would like to report on this conjecture in future publication.

\subsection*{Acknowledgements}
This work is supported by the Grant-in-Aid for 
Scientific Research {\bf C} (26400105)
from Japan Society for Promotion of Science.
The author would like to thank Professor Michio Jimbo for giving advice.
The author would like to thank Professor Zengo Tsuboi, Professor Pascal Baseilhac, Professor Kouichi Takemura and Professor Kenji Iohara for discussion.
The author would like to thank kind hospitality extended to him at University of Tours, University of Leeds and University of Lyon 1.

\begin{appendix}

\section{Normal ordering rules}

\label{appendix:A}

In this Appendix we summarize normal ordering rules that are useful for proof of main results.

\subsection{$[X^{+,i}(z_1),X^{-,i+1}(z_2)]~~~(1\leq i \leq M+N-2)$}
$\bullet$~For $1\leq i \leq M-1$ we have
\begin{eqnarray}
[E_{i,j}^\pm(z_1),F_{i+1,j}^{1,\pm}(z_2)]=\pm (q-q^{-1})\delta(q^{-k+1-2j}z_2/z_1)
:E_{i,j}^\pm(z_1)F_{i+1,j}^{1,\pm}(z_2):~~(1\leq j \leq i).
\end{eqnarray}
\begin{eqnarray}
:E_{i,j}^+(q^{-k+1-2j}z)F_{i+1,j}^{1,+}(z):=:E_{i,j}^-(q^{-k+1-2j}z)F_{i+1,j}^{1,-}(z):~~~(1\leq j \leq i).
\end{eqnarray}
$\bullet$~For $i=M$ we have
\begin{eqnarray}
[E_{M,j}(z_1),F_{M+1,j}^1(z_2)]=0~~~(1\leq j \leq M).
\end{eqnarray}
$\bullet$~For $M+1 \leq i \leq M+N-2$ we have
\begin{eqnarray}
&&[E_{i,j}(z_1),F_{i+1,j}^1(z_2)]=0~~~(1\leq j \leq M),
\\
&&[E_{i,j}^\pm(z_1),F_{i+1,j}^{1,\pm}(z_2)]=\mp (q-q^{-1})\delta(q^{-k-4M-1+2j}z_2/z_1)
:E_{i,j}^\pm(z_1)F_{i+1,j}^{1,\pm}(z_2):~(M+1 \leq j \leq i).\nonumber
\\
\\
&&:E_{i,j}^+(q^{-k-4M-1+2j}z)F_{i+1,j}^{1,+}(z):
=:E_{i,j}^-(q^{-k-4M-1+2j}z)F_{i+1,j}^{1,-}(z):~~~(M+1 \leq j \leq i).
\end{eqnarray}

\subsection{$[X^{+,i}(z_1),X^{-,i-1}(z_2)]~~~(2\leq i \leq M+N-1)$}
$\bullet$~For $2\leq i \leq M-1$ we have
\begin{eqnarray}
&&[E_{i,i-1}^+(z_1),F_{i-1,i+1}^{3,\pm}(z_2)]=-(q-q^{-1})\delta(q^{k+1}z_2/z_1)
:E_{i,i-1}^+(z_1)F_{i-1,i+1}^{3,\pm}(z_2):,\\
&&[E_{i,i}^\pm(z_1),F_{i-1,i-1}^{2,+}(z_2)]=(q-q^{-1})\delta(q^{k+1}z_2/z_1)
:E_{i,i}^\pm(z_1)F_{i-1,i-1}^{2,+}(z_2):,\\
&&[E_{i,i}^\pm(z_1),F_{i-1,i+1}^{3,\mp}(z_2)]=\pm(q-q^{-1})\delta(q^{k+1}z_2/z_1)
:E_{i,i}^\pm(z_1)F_{i-1,i+1}^{3,\mp}(z_2):.
\end{eqnarray}
\begin{eqnarray}
&&:E_{i,i-1}^+(q^{k+1}z)F_{i-1,i+1}^{3,\pm}(z):=:E_{i,i}^\pm(q^{k+1}z)F_{i-1,i-1}^{2,+}(z):,\\
&&:E_{i,i}^+(q^{k+1}z)F_{i-1,i+1}^{3,-}(z):=
:E_{i,i}^-(q^{k+1}z)F_{i-1,i+1}^{3,+}(z):.
\end{eqnarray}
$\bullet$~For $i=M$ we have
\begin{eqnarray}
&&[E_{M,M-1}(z_1),F_{M-1,M+1}^3(z_2)]=\frac{1}{q^{M-2}z_1}\delta(q^{k+1}z_2/z_1)
:E_{M,M-1}(z_1)F_{M-1,M+1}^3(z_2):,\\
&&[E_{M,M}(z_1),F_{M-1,M-1}^{2,+}(z_2)]=(q-q^{-1})\delta(q^{k+1}z_2/z_1)
:E_{M,M}(z_1)F_{M-1,M-1}^{2,+}(z_2):,\\
&&[E_{M,M}(z_1),F_{M-1,M+1}^3(z_2)]=0.
\end{eqnarray}
\begin{eqnarray}
:E_{M,M-1}(q^{k+1}z)F_{M-1,M+1}^3(z):=:E_{M,M}(q^{k+1}z)F_{M-1,M-1}^{2,+}(z):.
\end{eqnarray}
$\bullet$~For $i=M+1$ we have
\begin{eqnarray}
&&[E_{M+1,M}(z_1),F_{M,M+2}^{3,\pm}(z_2)]=\frac{1}{q^{M-1}z_1}\delta(q^{k-1}z_2/z_1)
:E_{M+1,M}(z_1)F_{M,M+2}^{3,\pm}(z_2):,\\
&&[E_{M+1,M+1}^\pm(z_1),F_{M,M}^{2,+}(z_2)]=-(q-q^{-1})\delta(q^{k-1}z_2/z_1)
:E_{M+1,M+1}^\pm(z_1)F_{M,M}^{2,+}(z_2):,
\\
&&[E_{M+1,M+1}^\pm(z_1),F_{M,M+2}^{3,\mp}(z_2)]=\mp q(q-q^{-1})\delta(q^{k-1}z_2/z_1)
:E_{M+1,M+1}^\pm(z_1)F_{M,M+2}^{3,\mp}(z_2):.
\end{eqnarray}
\begin{eqnarray}
&&:E_{M+1,M}(q^{k-1}z)F_{M,M+2}^{3,\pm}(z):=:E_{M+1,M+1}^\pm(q^{k-1}z)F_{M,M}^{2,+}(z):,\\
&&:E_{M+1,M+1}^+(q^{k-1}z)F_{M,M+2}^{3,-}(z):=
:E_{M+1,M+1}^-(q^{k-1}z)F_{M,M+2}^{3,+}(z):.
\end{eqnarray}
$\bullet$~For $M+2\leq i \leq M+N-1$ we have
\begin{eqnarray}
&&[E_{i,i-1}^+(z_1),F_{i-1,i+1}^{3,\pm}(z_2)]=(q-q^{-1})\delta(q^{k-1}z_2/z_1):
E_{i,i-1}^+(z_1)F_{i-1,i+1}^{3,\pm}(z_2):,\\
&&[E_{i,i}^\pm(z_1),F_{i-1,i-1}^{2,+}(z_2)]=-(q-q^{-1})\delta(q^{k-1}z_2/z_1):
E_{i,i}^\pm(z_1)F_{i-1,i-1}^{2,+}(z_2):,
\\
&&[E_{i,i}^\pm(z_1),F_{i-1,i+1}^{3,\mp}(z_2)]=\mp q(q-q^{-1})\delta(q^{k-1}z_2/z_1):
E_{i,i}^\pm(z_1)F_{i-1,i+1}^{3,\mp}(z_2):.
\end{eqnarray}
\begin{eqnarray}
&&:E_{i,i-1}^+(q^{k-1}z)F_{i-1,i+1}^{3,\pm}(z):=:E_{i,i}^\pm(q^{k-1}z)F_{i-1,i-1}^{2,+}(z):,\\
&&:E_{i,i}^+(q^{k-1}z)F_{i-1,i+1}^{3,-}(z):=:E_{i,i}^-(q^{k-1}z)F_{i-1,i+1}^{3,+}(z):.
\end{eqnarray}

\subsection{$[X^{+,i}(z_1),X^{-,j}(z_2)]~~~(3\leq i \leq M+N-1, 1\leq j \leq i-2)$}
$\bullet$~For $3\leq i \leq M-1$ and $1\leq j \leq i-2$ we have
\begin{eqnarray}
&&[E_{i,j}^+(z_1),F_{j,i+1}^{3,\pm}(z_2)]=-(q-q^{-1})\delta(q^{k+i-j}z_2/z_1):
E_{i,j}^+(z_1)F_{j,i+1}^{3,\pm}(z_2):,
\\
&&[E_{i,j+1}^\pm(z_1),F_{j,i}^{3,+}(z_2)]=(q-q^{-1})\delta(q^{k+i-j}z_2/z_1)
:E_{i,j+1}^\pm(z_1)F_{j,i}^{3,+}(z_2):,
\\
&&[E_{i,j+1}^\pm(z_1),F_{j,i+1}^{3,\mp}(z_2)]=\pm (q-q^{-1})\delta(q^{k+i-j}z_2/z_1):
E_{i,j+1}^\pm(z_1)F_{j,i+1}^{3,\mp}(z_2):.
\end{eqnarray}
\begin{eqnarray}
&&:E_{i,j}^+(q^{k+i-j}z)F_{j,i+1}^{3,\pm}(z):=:E_{i,j+1}^\pm(q^{k+i-j}z)F_{j,i}^{3,+}(z):,\\
&&:E_{i,j+1}^+(q^{k+i-j}z)F_{j,i+1}^{3,-}(z):=:E_{i,j-1}^-(q^{k+i-j}z)F_{j,i+1}^{3,+}(z):.
\end{eqnarray}
$\bullet$~For $i=M$ and $1\leq j \leq M-2$ we have
\begin{eqnarray}
&&[E_{M,j+1}(z_1),F_{j,M}^{3,+}(z_2)]=(q-q^{-1})\delta(q^{k+M-j}z_2/z_1):
E_{M,j+1}(z_1)F_{j,M}^{3,+}(z_2):,
\\
&&[E_{M,j}(z_1),F_{j,M+1}^3(z_2)]=\frac{1}{q^{j-1}z}\delta(q^{k+M-j}z_2/z_1):
E_{M,j}(z_1)F_{j,M+1}^3(z_2):,
\\
&&[E_{M,j+1}(z_1),F_{j,M+1}^3(z_2)]=0.
\end{eqnarray}
\begin{eqnarray}
:E_{M,j+1}(q^{k+M-j}z)F_{j,M}^{3,+}(z):=
:E_{M,j}(q^{k+M-j}z)F_{j,M+1}^3(z):.
\end{eqnarray}
$\bullet$~For $M+1 \leq i \leq M+N-1$ and $1\leq j \leq M-1$ we have
\begin{eqnarray}
&&[E_{i,j+1}(z_1),F_{j,i}^3(z_2)]=-\frac{1}{q^{j+1}z_1}\delta(q^{k+2M-i-j}z_2/z_1):
E_{i,j+1}(z_1)F_{j,i}^3(z_2):,
\\
&&[E_{i,j}(z_1),F_{j,i+1}^3(z_2)]=\frac{1}{q^{j-1}z_1}\delta(q^{k+2M-i-j}z_2/z_1):
E_{i,j}(z_1)F_{j,i+1}^3(z_2):,\\
&&[E_{i,j+1}(z_1),F_{j,i+1}^3(z_2)]=0.
\end{eqnarray}
\begin{eqnarray}
:E_{i,j+1}(q^{k+2M-i-j}z)F_{j,i}^3(z):=
:E_{i,j}(q^{k+2M-i-j}z)F_{j,i+1}^3(z):.
\end{eqnarray}
$\bullet$~For $M+2 \leq i \leq M+N-1$ and $j=M$ we have
\begin{eqnarray}
&&[E_{i,M+1}^\pm(z_1),F_{M,i}^{3,+}(z_2)]=-(q-q^{-1})\delta(q^{k+M-i}z_2/z_1):
E_{i,M+1}^\pm(z_1)F_{M,i}^{3,+}(z_2):,
\\
&&[E_{i,M}(z_1),F_{M,i+1}^{3,\pm}(z_2)]=\frac{1}{q^{M-1}z_1}\delta(q^{k+M-i}z_2/z_1):
E_{i,M}(z_1)F_{M,i+1}^{3,\pm}(z_2):,
\\
&&[E_{i,M+1}^\pm(z_1),F_{M,i+1}^{3,\mp}(z_2)]=\mp q(q-q^{-1})\delta(q^{k+M-i}z_2/z_1):
E_{i,M+1}^\pm(z_1)F_{M,i+1}^{3,\mp}(z_2):.
\end{eqnarray}
\begin{eqnarray}
&&
:E_{i,M+1}^\pm(q^{k+M-i}z)F_{M,i}^{3,+}(z):
=:E_{i,M}(q^{k+M-i}z)F_{M,i+1}^{3,\pm}(z):,\\
&&:E_{i,M+1}^+(q^{k+M-i}z)F_{M,i+1}^{3,-}(z):=
:E_{i,M-1}^-(q^{k+M-i}z)F_{M,i+1}^{3,-}(z):.
\end{eqnarray}
$\bullet$~For $M+2 \leq i \leq M+N-1$ and $M+1 \leq j \leq i-2$ we have
\begin{eqnarray}
&&[E_{i,j+1}^{\pm}(z_1),F_{j,i}^{3,+}(z_2)]=-(q-q^{-1})\delta(q^{k-i+j}z_2/z_1):
E_{i,j+1}^{\pm}(z_1)F_{j,i}^{3,+}(z_2):,
\\
&&[E_{i,j}^+(z_1),F_{j,i+1}^{3,\pm}(z_2)]=(q-q^{-1})\delta(q^{k-i+j}z_2/z_1):
E_{i,j}^+(z_1)F_{j,i+1}^{3,\pm}(z_2):,
\\
&&[E_{i,j+1}^\pm(z_1),F_{j,i+1}^{3,\mp}(z_2)]=\mp q(q-q^{-1})\delta(q^{k-i+j}z_2/z_1):
E_{i,j+1}^\pm(z_1)F_{j,i+1}^{3,\mp}(z_2):.
\end{eqnarray}
\begin{eqnarray}
&&
:E_{i,j+1}^{\pm}(q^{k-i+j}z)F_{j,i}^{3,+}(z):=:E_{i,j}^+(q^{k-i+j}z)F_{j,i+1}^{3,\pm}(z):,
\\
&&:E_{i,j+1}^+(q^{k-i+j}z)F_{j,i+1}^{3,-}(z):=:E_{i,j+1}^-(q^{k-i+j}z)F_{j,i+1}^{3,+}(z):.
\end{eqnarray}

\subsection{$[X^{-,i}(z_1),X^{-,i-1}(z_2)]~~~(2\leq i \leq M+N-1)$}
$\bullet$~For $2 \leq i \leq M-1$ we have
\begin{eqnarray}
&&[F_{i,i-1}^{1,\pm}(z_1),F_{i-1,i-1}^{2,\pm}(z_2)]=\pm (q-q^{-1})\delta(q^{2k+2i-1}z_2/z_1):
F_{i,i-1}^{1,\pm}(z_1)F_{i-1,i-1}^{2,\pm}(z_2):.
\end{eqnarray}
\begin{eqnarray}
:F_{i,i-1}^{1,+}(q^{2k+2i-1}z)F_{i-1,i-1}^{2,+}(z):
=:F_{i,i-1}^{1,-}(q^{2k+2i-1}z)F_{i-1,i-1}^{2,-}(z):.
\end{eqnarray}
$\bullet$~For $i=M$ we have
\begin{eqnarray}
&&[F_{M,M-1}^{1,+}(z_1),F_{M-1,M-1}^{2,+}(z_2)]=(q-q^{-1})\delta(q^{2k+2M-1}z_2/z_1):
F_{M,M-1}^{1,+}(z_1)F_{M-1,M-1}^{2,+}(z_2)
:,
\\
&&[F_{M,M}^{2,-}(z_1),F_{M-1,M+1}^3(z_2)]=-\frac{q^{k+M-1}}{z_1}\delta(q^{2k+2M-1}z_2/z_1):
F_{M,M}^{2,-}(z_1)F_{M-1,M+1}^3(z_2):.
\end{eqnarray}
\begin{eqnarray}
:F_{M,M-1}^{1,+}(q^{2k+2M-1}z)F_{M-1,M-1}^{2,+}(z):=
:F_{M,M}^{2,-}(q^{2k+2M-1}z)F_{M-1,M+1}^3(z):.
\end{eqnarray}
$\bullet$~For $i=M+1$ we have
\begin{eqnarray}
&&[F_{M+1,M}^1(z_1),F_{M,M}^{2,+}(z_2)]=-\frac{q^{k+M}}{z_1}\delta(q^{2k+2M-1}z_2/z_1)
:F_{M+1,M}^1(z_1)F_{M,M}^{2,+}(z_2):,
\\
&&[F_{M+1,M+1}^{2,-}(z_1),F_{M,M+2}^{3,-}(z_2)]=(q-q^{-1})\delta(q^{2k+2M-1}z_2/z_1):
F_{M+1,M+1}^{2,-}(z_1)F_{M,M+2}^{3,-}(z_2):.
\end{eqnarray}
\begin{eqnarray}
:F_{M+1,M}^1(q^{2k+2M-1}z)F_{M,M}^{2,+}(z):=
:F_{M+1,M+1}^{2,-}(q^{2k+2M-1}z)F_{M,M+2}^{3,-}(z):.
\end{eqnarray}
$\bullet$~For $M+2 \leq i \leq M+N-1$ we have
\begin{eqnarray}
&&[F_{i,i-1}^{1,+}(z_1),F_{i-1,i-1}^{2,+}(z_2)]=\mp
(q-q^{-1})\delta(q^{2k+4M-2i+1}z_2/z_1):
F_{i,i-1}^{1,+}(z_1)F_{i-1,i-1}^{2,+}(z_2):.
\end{eqnarray}
\begin{eqnarray}
:F_{i,i-1}^{1,+}(q^{2k+4M-2i+1}z)F_{i-1,i-1}^{2,+}(z):
=:F_{i,i-1}^{1,-}(q^{2k+4M-2i+1}z)F_{i-1,i-1}^{2,-}(z):.
\end{eqnarray}

\subsection{$[X^{-,i}(z_1),X^{-,j}(z_2)]~~~(3\leq i \leq M+N-1, 1\leq j \leq i-2)$}
$\bullet$ For $3\leq i \leq M-1$ and $1\leq j \leq i-2$ we have
\begin{eqnarray}
&&[F_{i,j}^{1,+}(z_1),F_{j,i}^{3,\pm}(z_2)]=q(q-q^{-1})\delta(q^{2k+i+j}z_2/z_1):
F_{i,j}^{1,+}(z_1)F_{j,i}^{3,\pm}(z_2):,
\\
&&[F_{i,j+1}^{1,\pm}(z_1),F_{j,i+1}^{3,-}(z_2)]
=-(q-q^{-1})\delta(q^{2k+i+j}z_2/z_1):
F_{i,j+1}^{1,\pm}(z_1)F_{j,i+1}^{3,-}(z_2):,
\\
&&[F_{i,j+1}^{1,\pm}(z_1),F_{j,i}^{3,\pm}(z_2)]
=\mp (q-q^{-1})\delta(q^{2k+i+j}z_2/z_1):
F_{i,j+1}^{1,\pm}(z_1)F_{j,i}^{3,\pm}(z_2):.
\end{eqnarray}
\begin{eqnarray}
&&:F_{i,j}^{1,+}(q^{2k+i+j}z)F_{j,i}^{3,\pm}(z):
=:F_{i,j+1}^{1,\pm}(q^{2k+i+j}z)F_{j,i+1}^{3,-}(z):,
\\
&&:F_{i,j+1}^{1,+}(q^{2k+i+j}z)F_{j,i}^{3,+}(z):
=:F_{i,j+1}^{1,-}(q^{2k+i+j}z)F_{j,i}^{3,-}(z):.
\end{eqnarray}
$\bullet$~For $i=M$ and $1\leq j \leq M-2$ we have
\begin{eqnarray}
&&[F_{M,j}^{1,+}(z_1),F_{j,M}^{3,\pm}(z_2)]=q(q-q^{-1})\delta(q^{2k+M+j}z_2/z_1):
F_{M,j}^{1,+}(z_1)F_{j,M}^{3,\pm}(z_2):,
\\
&&[F_{M,j+1}^{1,\pm}(z_1),F_{j,M+1}^3(z_2)]
=-\frac{q^{k+j-1}}{z_1}\delta(q^{2k+M+j}z_2/z_1):
F_{M,j+1}^{1,\pm}(z_1)F_{j,M+1}^3(z_2):,
\\
&&[F_{M,j+1}^{1,\pm}(z_1),F_{j,M}^{3,\pm}(z_2)]=\mp (q-q^{-1})\delta(q^{2k+M+j}z_2/z_1):
F_{M,j+1}^{1,\pm}(z_1)F_{j,M}^{3,\pm}(z_2):.
\end{eqnarray}
\begin{eqnarray}
&&:F_{M,j}^{1,+}(q^{2k+M+j}z)F_{j,M}^{3,\pm}(z):=:F_{M,j+1}^{1,\pm}(q^{2k+M+j}z)F_{j,M+1}^{3}(z):,
\\
&&:F_{M,j+1}^{1,+}(q^{2k+M+j}z)F_{j,M}^{3,+}(z):
=:F_{M,j+1}^{1,-}(q^{2k+M+j}z)F_{j,M}^{3,-}(z):.
\end{eqnarray}
$\bullet$~For $M+1 \leq i \leq M+N-1$ and $M+1 \leq j \leq i-2$ we have
\begin{eqnarray}
&&[F_{i,j}^{1,+}(z_1),F_{j,i}^{3,\pm}(z_2)]
=-q^{-1}(q-q^{-1})\delta(q^{2k+4M-i-j}z_2/z_1):
F_{i,j}^{1,+}(z_1)F_{j,i}^{3,\pm}(z_2):,
\\
&&[F_{i,j+1}^{1,\pm}(z_1),F_{j,i+1}^{3,-}(z_2)]=(q-q^{-1})\delta(q^{2k+4M-i-j}z_2/z_1)
:F_{i,j+1}^{1,\pm}(z_1)F_{j,i+1}^{3,-}(z_2):,
\\
&&[F_{i,j+1}^{1,\pm}(z_1),F_{j,i}^{3,\pm}(z_2)]=\pm (q-q^{-1})\delta(q^{2k+4M-i-j}z_2/z_1)
:F_{i,j+1}^{1,\pm}(z_1)F_{j,i}^{3,\pm}(z_2):.
\end{eqnarray}
\begin{eqnarray}
&&:F_{i,j}^{1,+}(q^{2k+4M-i-j}z)F_{j,i}^{3,\pm}(z):=:F_{i,j+1}^{1,\pm}(q^{2k+4M-i-j}z)F_{j,i+1}^{3,-}(z):,
\\
&&:F_{i,j+1}^{1,+}(q^{2k+4M-i-j}z)F_{j,i}^{3,+}(z):=:F_{i,j+1}^{1,-}(q^{2k+4M-i-j}z)F_{j,i}^{3,-}(z):.
\end{eqnarray}
$\bullet$~For $M+2 \leq i \leq M+N-1$ and $j=M$ we have
\begin{eqnarray}
&&[F_{i,M}^1(z_1),F_{M,i}^{3,\pm}(z_2)]
=-\frac{q^{k+M-1}}{z_1}\delta(q^{2k+3M-i}z_2/z_1)
:F_{i,M}^1(z_1)F_{M,i}^{3,\pm}(z_2):,\\
&&[F_{i,M+1}^{1,\pm}(z_1),F_{M,i+1}^{3,-}(z_2)]=(q-q^{-1})\delta(q^{2k+3M-i}z_2/z_1):
F_{i,M+1}^{1,\pm}(z_1)F_{M,i+1}^{3,-}(z_2)
:,\\
&&[F_{i,M+1}^{1,\pm}(z_1),F_{M,i}^{3,\pm}(z_2)]=
\pm (q-q^{-1})\delta(q^{2k+3M-i}z_2/z_1):
F_{i,M+1}^{1,\pm}(z_1)F_{M,i}^{3,\pm}(z_2):.
\end{eqnarray}
\begin{eqnarray}
&&:F_{i,M}^{1,+}(q^{2k+3M-i}z)F_{M,i}^{3,\pm}(z):=:F_{i,M+1}^{1,\pm}(q^{2k+3M-i}z)F_{M,i+1}^{3,-}(z):,
\\
&&:F_{i,M+1}^{1,+}(q^{2k+3M-i}z)F_{M,i}^{3,+}(z):=:F_{i,M+1}^{1,-}(q^{2k+i+j}z)F_{M,i}^{3,-}(z):.
\end{eqnarray}
$\bullet$~For $M+1 \leq i \leq M+N-1$ and $1 \leq j \leq M-1$ we have
\begin{eqnarray}
&&
[F_{i,j}^1(z_1),F_{j,i}^{3}(z_2)]=-\frac{q^{k+j-1}}{z_1}\delta(q^{2k+2M-i+j}z_2/z_1):
F_{i,j}^1(z_1)F_{j,i}^{3}(z_2):,
\\
&&[F_{i,j+1}^{1}(z_1),F_{j,i+1}^{3}(z_2)]=
\frac{q^{k+j+1}}{z_1}\delta(q^{2k+2M-i+j}z_2/z_1):
F_{i,j+1}^{1}(z_1)F_{j,i+1}^{3}(z_2):,
\\
&&[F_{i,j+1}^{1}(z_1),F_{j,i}^{3}(z_2)]=0.
\end{eqnarray}
\begin{eqnarray}
:F_{i,j}^1(q^{2k+2M-i+j}z)F_{j,i}^{3}(z):=
:F_{i,j+1}^{1}(q^{2k+2M-i+j}z)F_{j,i+1}^{3}(z):.
\end{eqnarray}

\subsection{$[S_i(z_1),X^{+,j}(z_2)]~~~(1\leq i \leq M+N-2, i+1 \leq j \leq M+N-1)$}
$\bullet$~For $1\leq i \leq M-2$ and $i+1 \leq j \leq M-1$ we have
\begin{eqnarray}
&&
[S_{i,j}^+(z_1),E_{j,i}^+(z_2)]=q(q-q^{-1})\delta(q^{i+j-M+N}z_2/z_1):
S_{i,j}^+(z_1)E_{j,i}^+(z_2):,\\
&&
[S_{i,j+1}^-(z_1),E_{j,i+1}^-(z_2)]=-q(q-q^{-1})\delta(q^{i+j-M+N}z_2/z_1):
S_{i,j+1}^-(z_1)E_{j,i+1}^-(z_2):,\\
&&
[S_{i,j+1}^\pm(z_1),E_{j,i}^\pm(z_2)]=\mp (q-q^{-1})\delta(q^{i+j-M+N}z_2/z_1):
S_{i,j+1}^\pm(z_1)E_{j,i}^\pm(z_2):,\\
&&
[S_{i,j}^-(z_1),E_{j,i}^-(z_2)]=q(q-q^{-1})\delta(q^{i+j-M+N}z_2/z_1):
S_{i,j}^-(z_1)E_{j,i}^-(z_2):,\\
&&
[S_{i,j+1}^+(z_1),E_{j,i+1}^-(z_2)]=-q(q-q^{-1})\delta(q^{i+j-M+N}z_2/z_1):
S_{i,j+1}^+(z_1)E_{j,i+1}^-(z_2):.
\end{eqnarray}
\begin{eqnarray}
&&
:S_{i,j}^+(q^{i+j-M+N}z)E_{j,i}^+(z):=:S_{i,j+1}^-(q^{i+j-M+N}z)E_{j,i+1}^-(z):,\\
&&
:S_{i,j+1}^+(q^{i+j-M+N}z)E_{j,i}^+(z):=
:S_{i,j+1}^-(q^{i+j-M+N}z)E_{j,i}^-(z):,
\\
&&
:S_{i,j}^+(q^{i+j-M+N}z)E_{j,i}^-(z):=
:S_{i,j+1}^+(q^{i+j-M+N}z)E_{j,i+1}^-(z):.
\end{eqnarray}
$\bullet$~For $i=M-1$ and $M+1 \leq j \leq M+N-1$ we have
\begin{eqnarray}
&&
[S_{M-1,M}^+(z_1),E_{M,M-1}(z_2)]=q^{-1}(q-q^{-1})\delta(q^{M+N-1}z_2/z_1):
S_{M-1,M}^+(z_1)E_{M,M-1}(z_2):,
\\
&&
[S_{M-1,j}(z_1),E_{j,M-1}(z_2)]=-\frac{q^{M+N-j-1}}{z_1}\delta(q^{2M+N-j-1}z_2/z_1):
S_{M-1,j}(z_1)E_{j,M-1}(z_2):.
\end{eqnarray}
For $i=M-1$ and $M+1 \leq j \leq M+N-1$ we have
\begin{eqnarray}
&&
[S_{M-1,M}^\pm(z_1),E_{M-1,M-1}^\pm(z_2)]=\mp (q-q^{-1})\delta(q^{M+N-2}z_2/z_1):
S_{M-1,M}^\pm(z_1)E_{M-1,M-1}^\pm(z_2):,\nonumber\\
\\
&&
[S_{M-1,j}(z_1),E_{j,M-1}(z_2)]=-\frac{q^{M+N-j-1}}{z_1}\delta(q^{2M+N-j-1}z_2/z_1):
S_{M-1,j}(z_1)E_{j,M-1}(z_2):.
\end{eqnarray}
For $i=M-1$ and $M \leq j \leq M+N-1$ we have
\begin{eqnarray}
&&
[S_{M-1,j+1}(z_1),E_{j,M}(z_2)]=-\frac{q^{M+N-j-1}}{z_1}\delta(q^{2M+N-j-1}z_2/z_1):
S_{M-1,j+1}(z_1)E_{j,M}(z_2):.
\end{eqnarray}
\begin{eqnarray}
&&:S_{M-1,M}^+(q^{M+N-2}z)E_{M-1,M-1}^+(z):=:S_{M-1,M}^-(q^{M+N-2}z)E_{M-1,M-1}^-(z):,
\\
&&:S_{M-1,M}^+(q^{M+N-1}z)E_{M,M-1}(z):=
:S_{M-1,M+1}(q^{M+N-1}z)E_{M,M}(z):,\\
&&
:S_{M-1,j}(q^{2M+N-j-1}z)E_{j,M-1}(z):=
:S_{M-1,j+1}(q^{2M+N-j-1}z)E_{j,M}(z):\nonumber\\
&&~~~(M+1\leq j \leq M+N-1).
\end{eqnarray}
$\bullet$~For $i=M$ and $M+1\leq j \leq M+N-1$ we have
\begin{eqnarray}
&&
[S_{M,j}(z_1),E_{j,M}(z_2)]=\frac{q^{M+N-j-1}}{z_1}\delta(q^{2M+N-j}z_2/z_1):
S_{M,j}(z_1)E_{j,M}(z_2):,
\\
&&
[S_{M,j+1}(z_1),E_{j,M+1}^-(z_2)]=q^{-1}(q-q^{-1})\delta(q^{2M+N-j}z_2/z_1)
:S_{M,j+1}(z_1)E_{j,M+1}^-(z_2):,
\\
&&
[S_{M,j+1}(z_1),E_{j,M}(z_2)]=0.
\end{eqnarray}
\begin{eqnarray}
&&
:S_{M,j}(q^{2M+N-j}z)E_{j,M}(z):=:S_{M,j+1}(q^{2M+N-j}z)E_{j,M+1}^-(z):.
\end{eqnarray}
$\bullet$~For $M+1 \leq i \leq M+N-1$ and $i+1 \leq j \leq M+N-1$ we have
\begin{eqnarray}
&&
[S_{i,j}^+(z_1),E_{j,i}^+(z_2)]=-q^{-1}(q-q^{-1})\delta(q^{3M+N-i-j}z_2/z_1):
S_{i,j}^+(z_1)E_{j,i}^+(z_2):,\\
&&
[S_{i,j+1}^-(z_1),E_{j,i+1}^-(z_2)]=q^{-1}(q-q^{-1})\delta(q^{3M+N-i-j}z_2/z_1):
S_{i,j+1}^-(z_1)E_{j,i+1}^-(z_2):,\\
&&
[S_{i,j+1}^\pm(z_1),E_{j,i}^\pm(z_2)]=\pm (q-q^{-1})\delta(q^{3M+N-i-j}z_2/z_1):
S_{i,j+1}^\pm (z_1)E_{j,i}^\pm (z_2):,\\
&&
[S_{i,j}^+(z_1),E_{j,i}^-(z_2)]=-q^{-1}(q-q^{-1})\delta(q^{3M+N-i-j}z_2/z_1):
S_{i,j}^+(z_1)E_{j,i}^-(z_2):,\\
&&
[S_{i,j+1}^+(z_1),E_{j,i+1}^-(z_2)]=q^{-1}(q-q^{-1})\delta(q^{3M+N-i-j}z_2/z_1)
:S_{i,j+1}^+(z_1)E_{j,i+1}^-(z_2):.
\end{eqnarray}
\begin{eqnarray}
&&
:S_{i,j}^+(q^{3M+N-i-j}z)E_{j,i}^+(z):=:S_{i,j+1}^-(q^{3M+N-i-j}z)E_{j,i+1}^-(z):,
\\
&&
:S_{i,j+1}^+(q^{3M+N-i-j}z)E_{j,i}^+(z):=
:S_{i,j+1}^-(q^{3M+N-i-j}z)E_{j,i}^-(z):,
\\
&&
:S_{i,j}^+(q^{3M+N-i-j}z)E_{j,i}^-(z):=
:S_{i,j+1}^+(q^{3M+N-i-j}z)E_{j,i+1}^-(z):.
\end{eqnarray}

\subsection{$[S_i(z_1),X^{-,i-1}(z_2)]~~~(2\leq i \leq M+N-1)$}
$\bullet$~For $2\leq i \leq M-1$ and $i+1\leq j \leq M$ we have
\begin{eqnarray}
[S_{i,j}^\pm(z_1),F_{i-1,j}^{3,\pm}(z_2)]=\mp (q-q^{-1})\delta(q^{k+2j-M+N-1}z_2/z_1):
S_{i,j}^\pm(z_1)F_{i-1,j}^{3,\pm}(z_2):.
\end{eqnarray}
\begin{eqnarray}
:S_{i,j}^+(q^{k+2j-M+N-1}z)F_{i-1,j}^{3,+}(z):=
:S_{i,j}^-(q^{k+2j-M+N-1}z)F_{i-1,j}^{3,-}(z):.
\end{eqnarray}
$\bullet$~For $i=M$ and $M+1 \leq j \leq M+N$ we have
\begin{eqnarray}
[S_{M,j}(z_1),F_{M-1,j}^3(z_2)]=0.
\end{eqnarray}
$\bullet$~For $i=M+1$ and $M+2 \leq j \leq M+N$ we have
\begin{eqnarray}
[S_{M+1,j}^\pm(z_1),F_{M,j}^{3,\pm}(z_2)]=\pm (q-q^{-1})\delta(q^{k+3M+N-2j+1}z_2/z_1):
S_{M+1,j}^\pm(z_1)F_{M,j}^{3,\pm}(z_2):.
\end{eqnarray}
\begin{eqnarray}
:S_{M+1,j}^+(q^{k+3M+N-2j+1}z)F_{M,j}^{3,+}(z):
=:S_{M+1,j}^-(q^{k+3M+N-2j+1}z)F_{M,j}^{3,-}(z):.
\end{eqnarray}
$\bullet$~For $M+2 \leq i \leq M+N-1$ and $i+1 \leq j \leq M+N$ we have
\begin{eqnarray}
[S_{i,j}^\pm(z_1),F_{i-1,j}^{3,\pm}(z_2)]=\pm (q-q^{-1})\delta(q^{k+3M+N-2j+1}z_2/z_1):
S_{i,j}^\pm(z_1)F_{i-1,j}^{3,\pm}(z_2):.
\end{eqnarray}
\begin{eqnarray}
:S_{i,j}^+(q^{k+3M+N-2j+1}z)F_{i-1,j}^{3,+}(z):
=:S_{i,j}^-(q^{k+3M+N-2j+1}z)F_{i-1,j}^{3,-}(z):.
\end{eqnarray}

\subsection{$[S_i(z_1),X^{-,j}(z_2)]~~~(1\leq i \leq M+N-2, i+1 \leq j \leq M+N-1)$}
$\bullet$~
For $1\leq i \leq M-1$ and $i+1\leq j \leq M$ we have
\begin{eqnarray}
~[S_{i,j+1}^-(z_1),F_{j,i}^{1,\pm}(z_2)]&=&
-(q-q^{-1})\delta(q^{-k-M+N-i+j}z_2/z_1):
S_{i,j+1}^-(z_1)F_{j,i}^{1,\pm}(z_2):,
\\
~[S_{i,j}^\pm(z_1),F_{j,i+1}^{1,-}(z_2)]&=&
(q-q^{-1})\delta(q^{-k-M+N-i+j}z_2/z_1):
S_{i,j}^\pm(z_1)F_{j,i+1}^{1,-}(z_2):
\nonumber\\
&&~~~(1\leq i \leq M-2, i+2 \leq j \leq M),\\
~[S_{i,i+1}^\pm(z_1),F_{i+1,i+1}^{2,-}(z_2)]
&=&(q-q^{-1})
\delta(q^{-k-M+N+1}z_2/z_1):
S_{i,i+1}^\pm(z_1)F_{i+1,i+1}^{2,-}(z_2):,
\\
~[S_{i,j}^\pm(z_1),F_{j,i}^{1,\mp}(z_2)]
&=&
\mp q^{-1}(q-q^{-1})
\delta(q^{-k-M+N-i+j}z_2/z_1):
S_{i,j}^\pm(z_1)F_{j,i}^{1,\mp}(z_2):\nonumber\\
&&~~~(1\leq i \leq M-2, i+1 \leq j \leq M-1),\\
~[S_{M-1,M}^\pm(z_1),F_{M,M-1}^{1,\mp}(z_2)]&=&
\pm q^{-1}(q-q^{-1})\delta(q^{-k-M+N+1}z_2/z_1):
S_{M-1,M}^\pm(z_1)F_{M,M-1}^{1,\mp}(z_2):.\nonumber\\
\end{eqnarray}
\begin{eqnarray}
:S_{i,i+1}^\pm(q^{-k-M+N+1}z)F_{i+1,i+1}^{2,-}(z):
&=&
:S_{i,i+2}^-(q^{-k-M+N+1}z)F_{i+1,i}^{1,\pm}(z):,\\
:S_{i,j}^\pm(q^{-k-M+N-i+j}z)F_{j,i+1}^{1,-}(z):
&=&
:S_{i,j+1}^-(q^{-k-M+N-i+j}z)F_{j,i}^{1,\pm}(z):\nonumber\\
&&~~~(1\leq i \leq M-2, i+2 \leq j \leq M),\\
:S_{i,j}^+(q^{-k-M+N-i+j}z)F_{j,i}^{1,-}(z):
&=&
:S_{i,j}^-(q^{-k-M+N-i+j}z)F_{j,i}^{1,+}(z):\nonumber\\
&&~~~(1\leq i \leq M-2, i+1 \leq j \leq M-1),
\\
:S_{M-1,M}^+(q^{-k-M+N+1}z)F_{M,M-1}^{1,-}(z):
&=&
:S_{M-1,M}^-(q^{-k-M+N+1}z)F_{M,M-1}^{1,+}(z):.
\end{eqnarray}
For $1 \leq i \leq M-1$ and $M+1 \leq j \leq M+N-1$ we have
\begin{eqnarray}
&&[S_{i,j}(z_1),F_{j,i}^1(z_2)]=0,
\\
&&[S_{i,j+1}(z_1),F_{j,i}^1(z_2)]=-\frac{q^{M+N-j-1}}{z_1}\delta(q^{-k+M+N-i-j}z_2/z_1):
S_{i,j+1}(z_1)F_{j,i}^1(z_2):,
\\
&&[S_{i,j}(z_1),F_{j,i+1}^1(z_2)]=\frac{q^{M+N-j+1}}{z_1}\delta(q^{-k+M+M-i-j}z_2/z_1):
S_{i,j}(z_1)F_{j,i+1}^1(z_2):.
\end{eqnarray}
\begin{eqnarray}
:S_{i,j}(q^{-k+M+N-i-j}z)F_{j,i+1}^1(z):
=:S_{i,j+1}(q^{-k+M+N-i-j}z)F_{j,i}^1(z):.
\end{eqnarray}
$\bullet$~For $i=M$ and $M+1 \leq j \leq M+N-1$ we have
\begin{eqnarray}
&&
[S_{M,j}(z_1),F_{j,M}^1(z_2)]=0,\\
&&
[S_{M,j+1}(z_1),F_{j,M}^1(z_2)]=
-\frac{q^{M+N-j-1}}{z_1}\delta(q^{-k+N-j}z_2/z_1):
S_{M,j+1}(z_1)F_{j,M}^1(z_2):,
\\
&&
[S_{M,j}(z_1),F_{j,M+1}^{1,-}(z_2)]=
-(q-q^{-1})\delta(q^{-k+N-j}z_2/z_1)
:S_{M,j}(z_1)F_{j,M+1}^{1,-}(z_2):.
\end{eqnarray}
\begin{eqnarray}
:S_{M,j+1}(q^{-k+N-j}z)F_{j,M}^1(z):=:S_{M,j}(q^{-k+N-j}z)F_{j,M+1}^{1,-}(z):.
\end{eqnarray}
$\bullet$~For $M+1 \leq i \leq M+N-1$ and $i+1\leq j \leq M+N-1$ we have
\begin{eqnarray}
~[S_{i,j+1}^-(z_1),F_{j,i}^{1,\pm}(z_2)]
&=&
(q-q^{-1})\delta(q^{-k-M+N+i-j}z_2/z_1)
:S_{i,j+1}^-(z_1)F_{j,i}^{1,\pm}(z_2):,
\\
~[S_{i,j}^\pm(z_1),F_{j,i+1}^{1,-}(z_2)]&=&
-(q-q^{-1})\delta(q^{-k-M+N+i-j}z_2/z_1)
:S_{i,j}^\pm(z_1)F_{j,i+1}^{1,-}(z_2):\nonumber\\
&&~~~(i+2 \leq j \leq M+N-1),\\
~[S_{i,i+1}^\pm(z_1),F_{i+1,i+1}^{2,-}(z_2)]
&=&-(q-q^{-1})
\delta(q^{-k-M+N+1}z_2/z_1):
S_{i,i+1}^\pm(z_1)F_{i+1,i+1}^{2,-}(z_2):,
\\
~[S_{i,j}^\pm(z_1),F_{j,i}^{1,\mp}(z_2)]
&=&\pm q^{-1}(q-q^{-1})
\delta(q^{-k-M+N+i-j}z_2/z_1):
S_{i,j}^\pm(z_1)F_{j,i}^{1,\mp}(z_2):.
\end{eqnarray}
For $M+1 \leq i \leq M+N-1$ and $i+2\leq j \leq M+N-1$ we have
\begin{eqnarray}
:S_{i,i+1}^\pm(q^{-k-M+N-1}z)F_{i+1,i+1}^{2,-}(z):
&=&
:S_{i,i+2}^-(q^{-k-M+N-1}z)F_{i+1,i}^{1,\pm}(z):,\\
:S_{i,j}^\pm(q^{-k-M+N+i-j}z)F_{j,i+1}^{1,-}(z):
&=&
:S_{i,j+1}^-(q^{-k-M+N+i-j}z)F_{j,i}^{1,\pm}(z):\nonumber\\
&&~~~(i+2 \leq j \leq M+N-1),
\\
:S_{i,j}^+(q^{-k-M+N+i-j}z)F_{j,i}^{1,-}(z):&=&
:S_{i,j}^-(q^{-k-M+N+i-j}z)F_{j,i}^{1,+}(z):.
\end{eqnarray}

\section{Difference realization of $U_q(sl(M|N))$}

\label{appendix:B}

In this Appendix we recall a $q$-difference realization of $U_q(sl(M|N))$ \cite{AOS2}.
We introduce the coordinates 
$x_{i,j}$ $(1\leq i<j \leq M+N)$ by
\begin{eqnarray}
x_{i,j}=\left\{\begin{array}{cc}
z_{i,j}&(\nu_i \nu_j=+),\\
\theta_{i,j}&(\nu_i \nu_j=-),
\end{array}
\right.
\end{eqnarray}
where $z_{i,j}$ are complex variables and
$\theta_{i,j}$ are the Grassmann odd variables
that satisfy
$\theta_{i,j}\theta_{i,j}=0$ and
$\theta_{i,j}\theta_{i',j'}=-\theta_{i',j'}\theta_{i,j}$.
We set the differential operators
$\vartheta_{i,j}=x_{i,j}
\frac{\partial}{\partial x_{i,j}}$.
We fix parameters $\lambda_i \in {\bf C}$, 
$(1\leq i \leq M+N-1)$.
We set $q$-difference operators
$h_i, e_i, f_i$ $(1\leq i \leq M+N-1)$ as follows.
\begin{eqnarray}
&&
h_i=\lambda_i+
\sum_{j=1}^{i-1}(\nu_i \vartheta_{j,i}-\nu_{i+1}\vartheta_{j,i+1})-(\nu_i+\nu_{i+1})\vartheta_{i,i+1}
+\sum_{j=i+1}^{M+N}
(\nu_{i+1}\vartheta_{i+1,j}-\nu_i \vartheta_{i,j}),
\\
&&
e_i=e_{i,i}+\sum_{j=1}^{i-1} e_{i,j},~~~
f_i=\nu_i \sum_{j=1}^{i-1} f_{i,j}^1+f_{i,i}^2-\nu_{i+1}\sum_{j=i+2}^{M+N}f_{i,j}^3,
\end{eqnarray}
where we have set
\begin{eqnarray}
e_{i,i}&=&\frac{1}{x_{i,i+1}}[\vartheta_{i,i+1}]_q
~q^{
\sum_{l=1}^{i-1}(\nu_i \vartheta_{l,i}-\nu_{i+1}
\vartheta_{l,i+1})},
\\
e_{i,j}&=&
x_{j,i}\frac{1}{x_{j,i+1}}[\vartheta_{j,i+1}]_q
~q^{
\sum_{l=1}^{j-1}(\nu_i \vartheta_{l,i}-\nu_{i+1}
\vartheta_{l,i+1})},\\
f_{i,j}^1&=&
x_{j,i+1}
\frac{1}{x_{j,i}}
[\vartheta_{j,i}]_q
q^{
\sum_{l=j+1}^{i-1}(\nu_{i+1}\vartheta_{l,i+1}-\nu_i \vartheta_{l,i})
-\lambda_i+
(\nu_i+\nu_{i+1})\vartheta_{i,i+1}
+\sum_{l=i+2}^{N+1}(\nu_i \vartheta_{i,l}-\nu_{i+1}\vartheta_{i+1,l})},
\\
f_{i,i}^2
&=&
x_{i,i+1}\left[\lambda_i-\nu_i \vartheta_{i,i+1}-
\sum_{l=i+2}^{N+1}
(\nu_i \vartheta_{i,l}-\nu_{i+1} \vartheta_{i+1,l})
\right]_q,\\
f_{i,j}^3&=&
x_{i,j+1}\frac{1}{x_{i+1,j}}
[\vartheta_{i+1,j}]_q q^{
\lambda_i+\sum_{l=j+1}^{N+1}
(\nu_{i+1}\vartheta_{i+1,l}-\nu_i \vartheta_{i,l})}.
\end{eqnarray}
For Grassmann odd variables
$x_{i,j}=\theta_{i,j}$, the expression $\frac{1}{x_{i,j}}$
stands for derivative $\frac{1}{x_{i,j}}=\frac{\partial}{\partial x_{i,j}}$.
The $q$-difference operators
$h_i, e_i, f_i$ $(1\leq i \leq M+N-1)$ satisfy
the defining relations of $U_q(sl(M|N))$.

\section{Limit $q \to 1$}
\label{appendix:C}

In this Appendix we summarize useful formulae to take the limit $q \to 1$. 
Using the following relations
\begin{eqnarray}
&&
b_\pm^{i,j}(z)-(b+c)^{i,j}(q^{\pm 1}z)=-c^{i,j}(q^{\pm 1}z)-b^{i,j}(z;-1),\\
&&
b_\pm^{i,j}(z)-(b+c)^{i,j}(q^{\mp 1}z)=-c^{i,j}(q^{\mp 1}z)-b^{i,j}(z;1),\\
&&
b_\pm^{i,j}(q^{\pm \alpha}z)=\left(\frac{1}{M}b^{i,j}\right)(q^{\pm M}z;\alpha)-\left(\frac{1}{M}b^{i,j}\right)(z;\alpha-M),\\
&&
b^{i,j}(q^{\pm \alpha}z)=\left(\frac{\alpha}{M}b^{i,j}\right)(q^{\pm M}z;0)+\left(\frac{M-\alpha}{M}b^{i,j}\right)(z;0),
\end{eqnarray}
we have
\begin{eqnarray}
&&
\frac{\pm 1}{(q-q^{-1})z}:(e^{\pm b_+^{i,j}(z)-(b+c)^{i,j}(q^{\pm 1}z)}-e^{\pm b_-^{i,j}(z)-(b+c)^{i,j}(q^{\mp 1}z)}):=
:{_1 \partial_z}\left(e^{-c^{i,j}(z)}\right) e^{-b^{i,j}(z;\mp 1)}:,\\
&&
\frac{1}{(q-q^{-1})z}:(e^{a_+^i(q^{\frac{k+g}{2}}z)}-e^{a_-^i(q^{-\frac{k+g}{2}}z)}):
=:{_{k+g} \partial_z}\left(e^{
\left(\frac{1}{k+g}a^i\right)\left(z;\frac{k+g}{2}\right)
}\right)
e^{-\left(\frac{1}{k+g}a^i\right)\left(z;-\frac{k+g}{2}\right)}:,
\\
&&
\frac{1}{(q-q^{-1})z}:(e^{b^{i,j}(q^{k+l}z)}-e^{b^{i,j}(q^{-k-l}z)}):=:{_{k+g} \partial_z}\left(e^{\left(\frac{k+l}{k+g}b^{i,j}\right)(z;0)}\right)
e^{\left(\frac{g-l}{k+g}b^{i,j}\right)(z;0)}:,
\\
&&
\frac{1}{(q-q^{-1})z}:(e^{b_+^{i,j}(q^{k+l}z)}-e^{b_-^{i,j}(q^{-k-l}z)}):
=:{_{k+g} \partial_z}\left(e^{\left(\frac{1}{k+g}b^{i,j}\right)(z;k+l)}\right)
e^{-\left(\frac{1}{k+g}b^{i,j}\right)(z;l-g)}:.
\end{eqnarray}
Hence we have the following formulae.\\
$\bullet$
For $1\leq i \leq M-1$ we have
\begin{eqnarray}
H^i(z)&=&{_1 \partial_z}\left\{a^i\left(z;\frac{g}{2}\right)+\sum_{l=1}^i
\left(b^{l,i+1}\left(z;\frac{k}{2}+l-1\right)-b^{l,i}\left(z;\frac{k}{2}+l\right)\right)\right.\nonumber\\
&&\left.
-\sum_{l=i+1}^M
\left(b^{i+1,l}\left(z;\frac{k}{2}+l-1\right)-b^{i,l}\left(z;\frac{k}{2}+l\right)\right)\right.\\
&&\left.
-\sum_{l=M+1}^{M+N}\left(
b^{i+1,l}\left(z;\frac{k}{2}+2M-l\right)
-b^{i,l}\left(z;\frac{k}{2}+2M+1-l\right)
\right)
\right\}.\nonumber
\end{eqnarray}
$\bullet$
For $i=M$ we have
\begin{eqnarray}
H^M(z)&=&{_1 \partial_z}\left\{a^M\left(z;\frac{g}{2}\right)
-\sum_{l=1}^{M-1}
\left(b^{l,M+1}\left(z;\frac{k}{2}+l\right)+b^{l,M}\left(z;\frac{k}{2}+l\right)\right)\right.\nonumber\\
&&\left.
+\sum_{l=M+2}^{M+N}\left(
b^{M+1,l}\left(z;\frac{k}{2}+2M+1-l\right)
+b^{M,l}\left(z;\frac{k}{2}+2M+1-l\right)
\right)
\right\}.
\end{eqnarray}
$\bullet$
For $M+1\leq i \leq M+N-1$ we have 
\begin{eqnarray}
H^i(z)&=&{_1 \partial_z}\left\{a^i\left(z;\frac{g}{2}\right)
-\sum_{l=1}^M
\left(b^{l,i+1}\left(z;\frac{k}{2}+l\right)-b^{l,i}\left(z;\frac{k}{2}+l-1\right)\right)\right.\nonumber\\
&&\left.
-\sum_{l=M+1}^i
\left(b^{l,i+1}\left(z;\frac{k}{2}+2M-l+1\right)-b^{l,i}\left(z;\frac{k}{2}+2M-l\right)\right)\right.\\
&&\left.
+\sum_{l=i+1}^{M+N}\left(
b^{i+1,l}\left(z;\frac{k}{2}+2M-l+1\right)
-b^{i,l}\left(z;\frac{k}{2}+2M-l\right)
\right)
\right\}.\nonumber
\end{eqnarray}
$\bullet$
For $1\leq i \leq M-1$ we have
\begin{eqnarray}
&&
\frac{1}{(q-q^{-1})z}(E_{i,j}^+(z)-E_{i,j}^-(z))
\\
&=& :{_1 \partial_z}\left(e^{-c^{j,i+1}(q^{j-1}z)}\right) e^{-b^{j,i+1}(q^{j-1}z;-1)+(b+c)^{j,i}(q^{j-1}z)+\sum_{l=1}^{j-1}(\Delta_R^- b_+^{l,i})(q^lz)}:
~~~(1\leq j \leq i).\nonumber
\end{eqnarray}
$\bullet$
For $M+1 \leq i \leq M+N-1$ we have
\begin{eqnarray}
&&
\frac{1}{(q-q^{-1})z}(E_{i,j}^+(z)-E_{i,j}^-(z))\nonumber\\
&=&
-:{_1 \partial_z}\left(e^{-c^{j,i+1}(q^{2M+1-j}z)}\right)
e^{-b^{j,i+1}(q^{2M+1-j}z;1)+(b+c)^{j,i}(q^{2M+1-j}z)}\nonumber\\
&&\times
e^{-\sum_{l=1}^M(\Delta_R^+ b_+^{l,i}(q^{l-1}z))-\sum_{l=M+1}^{j-1}(\Delta_R^+ b_+^{l,i})(q^{2M-l}z)}: ~~~(M+1\leq j \leq i).
\end{eqnarray}
$\bullet$
For $1\leq i \leq M-1$ we have
\begin{eqnarray}
&&
\frac{-1}{(q-q^{-1})z}(F_{i,j}^{1,+}(z)-F_{i,j}^{1,-}(z))\nonumber\\
&=&
:{_1 \partial_z} \left(e^{-c^{j,i}(q^{-k-j}z)}\right) e^{-b^{j,i}(q^{-k-j}z;1)+a_-^i(q^{-\frac{k+g}{2}}z)+(b+c)^{j,i+1}(q^{-k-j}z)}\\
&&\times
e^{\sum_{l=j+1}^i (\Delta_R^+ b_-^{l,i})(q^{-k-l}z)-\sum_{l=i+1}^M(\Delta_L^+ b_-^{i,l})(q^{-k-l}z)-\sum_{l=M+1}^{M+N}
(\Delta_L^+ b_-^{i,l})(q^{-k-2M-1+l}z)}:~(1\leq j \leq i-1),\nonumber
\\
&&\frac{1}{(q-q^{-1})z}(F_{i,i}^{2,+}(z)-F_{i,i}^{2,-}(z))\nonumber\\
&=& 
:{_{k+g}\partial_z}
\left(e^{\left(\frac{1}{k+g}a^i\right)(z;\frac{k+g}{2})
+\left(\frac{k+i}{k+g}(b+c)^{i,i+1}\right)(z;0)
-\sum_{l=i+1}^M\left(
\left(\frac{1}{k+g}b^{i+1,l}\right)(z;k+l-1)
-
\left(\frac{1}{k+g}b^{i,l}\right)(z;k+l)
\right)
}\right.\nonumber\\
&&\times
\left.e^{-
\sum_{l=M+1}^{M+N}\left(
\left(\frac{1}{k+g}b^{i+1,l}\right)(z;k+2M-l)
-
\left(\frac{1}{k+g}b^{i,l}\right)(z;k+2M+1-l)
\right)
}\right)\nonumber\\
&&
\times
e^{
-\left(\frac{1}{k+g}a^i\right)(z;-\frac{k+g}{2})
+\left(\frac{g-i}{k+g}(b+c)^{i,i+1}\right)(z;0)
+\sum_{l=i+1}^M
\left(
\left(\frac{1}{k+g}b^{i+1,l}\right)(z;l-g-1)
-
\left(\frac{1}{k+g}b^{i,l}\right)(z;l-g)
\right)}\nonumber\\
&&
\times
e^{
\sum_{l=M+1}^{M+N}\left(
\left(\frac{1}{k+g}b^{i+1,l}\right)(z;2M-l-g)
-
\left(\frac{1}{k+g}b^{i,l}\right)(z;2M-l-g+1)
\right)}:,
\\
&&\frac{1}{(q-q^{-1})z}(F_{i,j}^{3,+}(z)-F_{i,j}^{3,-}(z))\nonumber\\
&=&
:{_1 \partial_z} \left(e^{-c^{i+1,j}(q^{k+j-1}z)}\right) e^{-b^{i+1,j}(q^{k+j-1}z;-1)
+a_+^i(q^{\frac{k+g}{2}}z)+(b+c)^{i,j}(q^{k+j-1}z)}\nonumber
\\
&&\times
e^{-\sum_{l=j}^M (\Delta_L^- b_+^{i,l})(q^{k+l}z)-\sum_{l=M+1}^{M+N}
(\Delta_L^- b_+^{i,l})(q^{k+2M+1-i}z)}:~~~(1\leq j \leq i-1).
\end{eqnarray}
$\bullet$
For $i=M$ we have
\begin{eqnarray}
&&
\frac{-1}{(q-q^{-1})z}(F_{M,j}^{1,+}(z)-F_{M,j}^{1,-}(z))
\nonumber\\
&=&
:{_1 \partial_z}\left(
e^{-c^{j,M}(q^{-k-j}z)}\right) e^{-b^{j,M}(q^{-k-j}z;1)
+a_-^M(q^{-\frac{k+g}{2}}z)
-b_-^{j,M+1}(q^{-k-j}z)
-b^{j,M+1}(q^{-k-j+1}z)}\\
&&\times
e^{-\sum_{l=j+1}^{M-1} (\Delta_R^0 b_-^{l,M})(q^{-k-l}z)
+\sum_{l=M+2}^{M+N}
(\Delta_L^0 b_-^{M,l})(q^{-k-2M-1+l}z)}:~~~(1\leq j \leq M-1),\nonumber
\\
&&\frac{1}{(q-q^{-1})z}(F_{M,M}^{2,+}(z)-F_{M,M}^{2,-}(z))\nonumber\\
&=& 
:{_{k+g}\partial_z}
\left(e^{\left(\frac{1}{k+g}a^M\right)(z;\frac{k+g}{2})
-\left(\frac{k+M-1}{k+g}b^{M,M+1}\right)(z;0)}\right.
\nonumber\\
&&\times \left.
e^{\sum_{l=M+2}^{M+N}
\left(
\left(\frac{1}{k+g}b^{M+1,l}\right)(z;k+2M+1-l)
+
\left(\frac{1}{k+g}b^{M,l}\right)(z;k+2M+1-l)
\right)
}\right.\nonumber\\
&&\times
\left.e^{-
\sum_{l=M+1}^{M+N}\left(
\left(\frac{1}{k+g}b^{M+1,l}\right)(z;k+2M-l)
-
\left(\frac{1}{k+g}b^{M,l}\right)(z;k+M+1-l)
\right)
}\right)\\
&&
\times
e^{
-\left(\frac{1}{k+g}a^M\right)(z;-\frac{k+g}{2})
-\left(\frac{g-M+1}{k+g}b^{M,M+1}\right)(z;0)
-\sum_{l=M+2}^{M+N}
\left(
\left(\frac{1}{k+g}b^{M+1,l}\right)(z;2M+1-l-g)
+
\left(\frac{1}{k+g}b^{M,l}\right)(z;2M+1-l-g)
\right)}:,
\nonumber
\\
&&
\frac{-1}{(q-q^{-1})z}(F_{M,j}^{3,+}(z)-F_{M,j}^{3,-}(z))\nonumber\\
&=&
:{_1 \partial_z}\left(e^{-c^{M+1,j}(q^{k+2M+1-j}z)}\right) 
e^{-b^{M+1,j}(q^{k+2M+1-j}z;1)
+a_+^M(q^{\frac{k+g}{2}}z)-b^{M,j}(q^{k+2M-j}z)}\nonumber
\\
&&\times
e^{b_+^{M+1,j}(q^{k+2M+1-j}z)+
\sum_{l=j+1}^{M+N} 
(\Delta_L^0 b_+^{M,l})(q^{k+2M+1-l}z)}:~~~(M+2 \leq j \leq M+N).
\end{eqnarray}
$\bullet$
For $M+1 \leq i \leq M+N-1$ we have
\begin{eqnarray}
&&
\frac{1}{(q-q^{-1})z}(F_{i,j}^{1,+}(z)-F_{i,j}^{1,-}(z))\nonumber\\
&=&
:{_1 \partial_z}\left(e^{-c^{j,i}(q^{-k-2M+j}z)}\right) e^{-b^{j,i}(q^{-k-2M+j}z;-1)+a_-^i(q^{-\frac{k+g}{2}}z)+(b+c)^{j,i+1}(q^{-k-2M+j}z)}\\
&&\times
e^{-\sum_{l=j+1}^i (\Delta_R^- b_-^{l,i})(q^{-k-2M+l}z)
+\sum_{l=i+1}^{M+N}(\Delta_L^- b_-^{i,l})(q^{-k-2M+l}z)}:~(M+1\leq j \leq i-1),\nonumber
\\
&&
\frac{1}{(q-q^{-1})z}(F_{i,i}^{2,+}(z)-F_{i,i}^{2,-}(z))\nonumber\\
&=& 
:{_{k+g}\partial_z}
\left(e^{\left(\frac{1}{k+g}a^i\right)(z;\frac{k+g}{2})
+\left(\frac{k+2M-i}{k+g}(b+c)^{i,i+1}\right)(z;0)
+
\sum_{l=i+1}^{M+N}\left(
\left(\frac{1}{k+g}b^{i+1,l}\right)(z;k+2M-l+1)
-
\left(\frac{1}{k+g}b^{i,l}\right)(z;k+2M-l)
\right)
}\right)\nonumber\\
&&
\times
e^{
-\left(\frac{1}{k+g}a^i\right)(z;-\frac{k+g}{2})
+\left(\frac{g-2M+i}{k+g}(b+c)^{i,i+1}\right)(z;0)
-\sum_{l=i+1}^{M+N}
\left(\left(\frac{1}{k+g}b^{i+1,l}\right)(z;2M-l-g+1)
-
\left(\frac{1}{k+g}b^{i,l}\right)(z;2M-l-g)
\right)},
\\
&&
\frac{-1}{(q-q^{-1})z}(F_{i,j}^{3,+}(z)-F_{i,j}^{3,-}(z))\nonumber\\
&=&
:{_1 \partial_z} \left(e^{-c^{i+1,j}(q^{k+2M-j+1}z)}\right) 
e^{-b^{i+1,j}(q^{k+2M-j+1}z;1)
+a_+^i(q^{\frac{k+g}{2}}z)+(b+c)^{i,j}(q^{k+2M-j+1}z)}\nonumber
\\
&&\times
e^{\sum_{l=j}^{M+N}
(\Delta_L^+ b_+^{i,l})(q^{k+2M-l}z)}:~~~(i+2 \leq j \leq M+N).
\end{eqnarray}
$\bullet$
For $1\leq i \leq M-1$ we have
\begin{eqnarray}
&&
\frac{-1}{(q-q^{-1})z}(\tilde{S}_{i,j}^+(z)-\tilde{S}_{i,j}^-(z))\nonumber\\
&=&
:{_1\partial_z}
\left(e^{-c^{i,j}(q^{M-N-j}z)}\right) 
e^{-b^{i,j}(q^{M-N-j}z;1)
+(b+c)^{i+1,j}(q^{M-N-j}z)}\nonumber
\\
&&\times
e^{\sum_{l=j+1}^M (\Delta_L^+ b_-^{i,l})(q^{M-N-l}z)
+\sum_{l=M+1}^{M+N}
(\Delta_L^+ b_-^{i,l})(q^{-M-N+l-1}z)}:
~~~
(i+1\leq j \leq M).
\end{eqnarray}
$\bullet$
For $M+1 \leq i \leq M+N-1$ we have
\begin{eqnarray}
&&
\frac{1}{(q-q^{-1})z}(\tilde{S}_{i,j}^+(z)-\tilde{S}_{i,j}^-(z))\nonumber\\
&=&
:{_1\partial_z}\left(e^{-c^{i,j}(q^{-M-N+j}z)}\right) e^{-b^{i,j}(q^{-M-N+j}z;-1)
+(b+c)^{i+1,j}(q^{-M-N+j}z)}\nonumber\\
&&\times e^{-\sum_{l=j+1}^{M+N} (\Delta_L^- b_-^{i,l})(q^{-M-N+l}z)}:
~~~
(i+1\leq j \leq M+N).
\end{eqnarray}
In the limit $q \to 1$
${_\alpha \partial_z}$, $\left(\frac{L}{M}a^{i}\right)(z;\alpha)$ and $\left(\frac{L}{M}b^{i,j}\right)(z;\alpha)$ become
$\alpha \partial_z$, $\frac{L}{M}a^i(z)$ and $\frac{L}{M} b^{i,j}(z)$ respectively.
Because the operators $a_\pm^i(z)$ $b_\pm^{i,j}(z)$, $(\Delta_L^\epsilon b_\pm^{i,j})(z)$ and
$(\Delta_R^\epsilon b_\pm^{i,j})(z)$ disappear,
our bosonization becomes simpler in the limit $q \to 1$.

\end{appendix}

\end{document}